\documentclass{amsart}

\usepackage{amsthm}
\usepackage{amssymb}
\usepackage{amsmath}
\usepackage{amsfonts}
\usepackage{amsrefs}
\usepackage{textcomp}
\usepackage[english]{babel}
\usepackage[cp1250]{inputenc}
\usepackage{textcomp}
\usepackage{wasysym}
\usepackage{stmaryrd}
\usepackage{esint}
\usepackage[all]{xy}
\usepackage{graphicx}

\newtheorem{tw}{Theorem}[section]
\newtheorem{dfn}[tw]{Definition}
\newtheorem{uw}[tw]{Remark}
\newtheorem{prz}[tw]{Example}
\newtheorem*{ozn}{Notation}
\newtheorem{lem}[tw]{Lemma}
\newtheorem{stw}[tw]{Proposition}
\newtheorem{wn}[tw]{Corollary}

\newtheorem*{dd}{Proof}
\newtheorem*{sd}{Idea of proof}
\newtheorem*{ak}{Acknowledgements}

\let\olddfn\dfn
\renewcommand{\dfn}{\olddfn\normalfont}

\let\oldlem\lem
\renewcommand{\lem}{\oldlem\normalfont}

\let\oldstw\stw
\renewcommand{\stw}{\oldstw\normalfont}

\let\oldozn\ozn
\renewcommand{\ozn}{\oldozn\normalfont}

\let\olduw\uw
\renewcommand{\uw}{\olduw\normalfont}

\let\oldwn\wn
\renewcommand{\wn}{\oldwn\normalfont}

\let\oldprz\prz
\renewcommand{\prz}{\oldprz\normalfont}

\let\olddd\dd
\renewcommand{\dd}{\olddd\normalfont}

\let\oldsd\sd
\renewcommand{\sd}{\oldsd\normalfont}

\let\oldak\ak
\renewcommand{\ak}{\oldak\normalfont}

\newcommand{\twopartdef}[4]
{
\left\{
		\begin{array}{ll}
			#1 & \mbox{if } #2 \\
			#3 & \mbox{if } #4
		\end{array}
	\right.
}

\newcommand{\threepartdef}[6]
{
	\left\{
		\begin{array}{lll}
			#1 & \mbox{if } #2 \\
			#3 & \mbox{if } #4 \\
			#5 & \mbox{if } #6
		\end{array}
	\right.
}

\author{Wojciech G\'{o}rny}

\address{W. G\'{o}rny: Faculty of Mathematics, Informatics and Mechanics, University of Warsaw, Warsaw, Poland.}

\email{w.gorny@mimuw.edu.pl}


\title{(Non)uniqueness of minimizers in the least gradient problem}


\begin{document}

\begin{abstract}
Minimizers in the least gradient problem with discontinuous boundary data need not be unique. However, all of them have a similar structure of level sets. Here, we give a full characterization of the set of minimizers in terms of any one of them and discuss stability properties of an approximate problem.
\end{abstract}

\maketitle

\section{Introduction}

Our main focus is the least gradient problem

\begin{equation}\label{zagadnienie}
\min \{ \int_\Omega |Du|, \quad u \in BV(\Omega), \quad Tu = f  \},
\end{equation}
where $T$ denotes the trace operator and $f \in L^1(\Omega)$. This paper deals with the issue of uniqueness of solutions to the least gradient problem. This type of problems, including anisotropic cases, has been adressed in many ways: from the point of view of geometric measure theory, see \cite{BGG}, \cite{SWZ}, \cite{JMN}, \cite{GRS}, via characterization of subdifferentials, see \cite{MRL}, \cite{Maz}, or as a reduction of a higher dimensional system coming from applications, namely conductivity imaging, again see \cite{JMN}, and free material design, again see \cite{GRS}.

In \cite{SWZ} it is estabilished that for continuous boundary data, under a condition on $\Omega$ slightly weaker than strict convexity, the solution exists and is continuous up to the boundary. Moreover, a maximum principle argument implies uniqueness of the minimizer. However, if we relax either continuity of boundary data or regularity properties of $\Omega$, we encounter additional difficulties:

(1) The solution itself might not exist: without strict convexity of $\Omega$ existence may fail even for continuous boundary data. This issue is discussed in \cite{GRS}, including some positive results on existence. On the other hand, as the example from \cite{ST} shows, if the boundary data belong only to $L^\infty(\partial\Omega)$, then the minimizer might not exist even if $\Omega$ is a two-dimensional disk. However, \cite{Gor} shows existence of solutions in the two-dimensional case for $BV$ boundary data.

(2) As pointed out in \cite{MRL}, uniqueness of solutions for discontinuous boundary data may fail even in the strictly convex case. However, all the solutions in their example have very similar structure of superlevel sets; they differ only on a set, on which each of the solutions is constant.

The paper is organized as follows: Section 2 provides the necessary background and some results concerning pointwise properties of precise representatives of least gradient functions. Section 3 is devoted to proving the main result of this paper, i.e. uniqueness of solutions to the least gradient problem except or a set where the solution is locally constant.

\begin{tw}\label{uniqueness}
Let $\Omega \subset \mathbb{R}^k$, where $2 \leq k \leq 7$, be an open bounded convex set with Lipschitz boundary. Let $u, v$ be precise representatives of functions of least gradient in $\Omega$ such that $Tu = Tv = h$. Then $u = v$ on $\Omega \backslash (C \cup N)$, where both $u$ and $v$ are locally constant on $C$ and $N$ has Hausdorff dimension at most $k-1$.
\end{tw}

Note that, unlike the existence results from \cite{SWZ} and \cite{Gor}, we require only convexity of $\Omega$ in place of some form of strict convexity of $\Omega$. However, we have an indirect assumption that the set $\Omega$ and the function $h$ support at least one solution to the least gradient problem. We do not address the question of necessary conditions for existence of solutions; an example of a set of sufficient conditions in $\mathbb{R}^2$, as given in \cite{Gor}, is that $\Omega$ is strictly convex with $C^1$ boundary and $h \in BV(\partial\Omega)$.

The proof will follow in two stages; firstly, the claim will be proved in the two-dimensional setting, where the proof faces less geometric difficulties. Then the claim will be proved for any $k$ such that the boundary of the superlevel set is an analytical minimal surface. This proof runs along similar lines, but with more serious geometrical difficulties and the two-dimensional proof will act as a toy model.

In Section 4 we use Theorem \ref{uniqueness} to provide a characterization of the set of solutions in terms of a single solution $u_0$. The results from this section are most useful in $\mathbb{R}^2$, as we consider certain partitions of sets by minimal surfaces; in dimensions higher than two finding all such partitions is a very hard question, while on the plane it can be turned into an algorithm.

Finally, Section 5 deals with an approximation of the least gradient problem which takes into account the total mass of the solution. Starting with $\Gamma-$convergence of corresponding functionals, we prove that minimizers of the approximate problems converge to a minimizers of least gradient problem with the smallest $L^p$ norm and this convergence is stronger that standard $L^p$ convergence.

\section{Preliminaries}

This section brings together a few technical results, which will be needed later, but are proved here not to interrupt the reasoning in section $3$. The general assumptions regarding the set $\Omega$ are the following: throughout the entire paper we will assume that $\Omega$ \textbf{is an open bounded set with Lipschitz boundary}. When necessary, we will impose the assumption of convexity of $\Omega$. Furthermore, in many results in Sections 2-4 we assume that $2 \leq k \leq 7$; this is necessary due to result by Giusti, see later in the commentary to Theorem \ref{twierdzeniezbgg}.

\subsection{Minimum of two BV functions}

The following two lemmas are simple exercises in BV theory. However, to the best of my knowledge, in the literature there is no proof for any of them. For more information regarding basic $BV$ theory, see \cite{AFP} or \cite{EG}. 

\begin{lem}\label{min+max}
Suppose that $u,v \in BV(\Omega)$. Then also $\min(u,v), \max(u,v) \in BV(\Omega)$ and the following inequality holds:

$$ \int_\Omega |D\max(u,v)| + \int_\Omega |D\min(u,v)| \leq \int_\Omega |Du| + \int_\Omega |Dv|.$$
\end{lem}

\begin{dd}
By \cite[Proposition 3.35]{AFP} we have for any sets $A, B$ of finite perimeter
$$ P(A \cup B, \Omega) + P(A \cap B, \Omega) \leq P(A, \Omega) + P(B, \Omega). $$
Let us plug into this inequality $A = E_t = \{u \geq t \}$ and $B = F_t = \{v \geq t\}$. Observe that $E_t \cup F_t = \{\max(u,v) \geq t\}$ and $E_t \cap F_t = \{\min(u,v) \geq t \}$. Thus for almost every $t$ (such that $E_t$ and $F_t$ have finite perimeter) we have
$$ P(\{\max(u,v) \geq t\}, \Omega) + P(\{\min(u,v) \geq t\}, \Omega) \leq P(\{u \geq t \}, \Omega) + P(\{u \geq t \}, \Omega). $$
Integration with respect to $t$ and the co-area formula give the result. \qed
\end{dd}

\begin{lem}\label{traceofminimum}
Let $u,v \in BV(\Omega)$. Then

$$ T\min(u,v) = \min(Tu, Tv) \qquad \mathcal{H}^{k-1}-a.e. \text{ on } \partial \Omega.$$
In particular, if $Tu = Tv = h$, then $T\min(u,v) = h$. Analogous result holds for $\max(u,v)$. 
\end{lem}

\begin{dd}
One inequality is obvious: the trace is a positive operator, so the inequality $\min(u,v) \leq u$ implies $T\min(u,v) \leq Tu$. Similarly $T\min(u,v) \leq Tv$, so $T\min(u,v) \leq \min(Tu, Tv)$.

For the opposite inequality, recall that for any $w \in BV(\Omega)$ on a set of full $\mathcal{H}^{k-1}$ measure we have

$$ \fint_{B(x,r) \cap \Omega} |w(y) - Tw(x)|dy \rightarrow 0.$$
Observe that this implies (by reverse triangle inequality for $L^1$ norm)

$$ 0 \leftarrow \fint_{B(x,r) \cap \Omega} |w(y) - Tw(x)|dy \geq |\fint_{B(x,r) \cap \Omega} |w(y)|dy - |Tw(x)||,$$
so 
$$\fint_{B(x,r) \cap \Omega} |w(y)|dy \rightarrow |Tw(x)|.$$
Now note that by linearity the trace of $w - s$ for $s \in \mathbb{R}$ equals $Tw - s$; thus for every $s \in \mathbb{R}$ we have

$$ \fint_{B(x,r) \cap \Omega} |w(y) - s|dy \rightarrow |Tw(x) - s|.$$
Let $P \subset \partial\Omega$ denote the set of full measure such that for every $x \in P$ the property above holds for $w = u, v, \min(u,v)$. Fix $x \in P$. Then we have

$$ \fint_{B(x,r) \cap \Omega} |u(y)|dy \rightarrow |a|, \, \fint_{B(x,r) \cap \Omega} |v(y)|dy \rightarrow |b|, \, \fint_{B(x,r) \cap \Omega} |\min(u,v)(y)|dy \rightarrow |c|.$$
Without loss of generality assume that $a \geq b$. We want to prove that $c \geq \min(a,b) = b$. We argue by contradiction: assume that $a \geq b > c$. Shift the functions $u, v$ by $s = a$, namely we obtain

$$ \fint_{B(x,r) \cap \Omega} |u(y) - a|dy \rightarrow |a - a| = 0, \, \fint_{B(x,r) \cap \Omega} |v(y) - a|dy \rightarrow |b - a|,$$
$$\fint_{B(x,r) \cap \Omega} |\min(u - a, v - a)(y)|dy \rightarrow |c - a|.$$
But $|\min(u - a, v - a)| \leq |u - a| + |v - a|$. Thus

$$|c - a| \leq \fint_{B(x,r) \cap \Omega} |\min(u - a, v - a)(y)|dy \leq \fint_{B(x,r) \cap \Omega} |u(y) - a|dy + $$
$$ \fint_{B(x,r) \cap \Omega} |v(y) - a|dy \rightarrow 0 + |b - a|,$$
but in the beginning we assumed that $a \geq b > c$, in particular $|c - a| > |b - a|$, contradiction. Thus $T\min(u,v)(x) \geq \min(Tu(x), Tv(x))$ for every $x \in P$, but it is a set of full measure. \qed
\end{dd}


\subsection{Least gradient functions}

In this subsection we recall the definition and some properties of least gradient functions; a standard reference is \cite{BGG} and \cite{Giu}. Then we prove some results concerning pointwise properties of precise representatives of least gradient functions.

\begin{dfn}
We say that $u \in BV(\Omega)$ is a function of least gradient, if for every compactly supported $($equivalently: with trace zero$)$ $v \in BV(\Omega)$ we have

\begin{equation*}
\int_\Omega |Du| \leq \int_\Omega |D(u + v)|.
\end{equation*}
We also say that $u$ is a solution of the least gradient problem for $f \in L^1(\partial\Omega)$ in the sense of traces, if $u$ is a least gradient function such that $Tu = f$.
\end{dfn}

To deal with regularity of least gradient functions, it is convenient to consider superlevel sets of $u$, i.e. sets of the form $\partial \{ u > t \}$ for $t \in \mathbb{R}$. A classical theorem states that

\begin{tw}\label{twierdzeniezbgg}
$($\cite[Theorem 1]{BGG}$)$ \\
Suppose $\Omega \subset \mathbb{R}^N$ is open. Let $u$ be a function of least gradient in $\Omega$. Then the set $\partial \{ u > t \}$ is minimal in $\Omega$, i.e. $\chi_{\{ u > t \}}$ is of least gradient for every $t \in \mathbb{R}$. \qed
\end{tw}
Obviously the theorem also holds for sets of the form $\{ u \geq t \}$. Similarly, all the results below could be stated for either $\{ u > t \}$ or $\{ u \geq t \}$; in Section $3$ we will use whichever version is more convenient.

This result was later improved in \cite[Chapter 10]{Giu} that in low dimensions $(k \leq 7)$ the boundary $\partial E$ of a minimal set $E$ is an analytical hypersurface $($after taking the precise representative of the set $E)$. In particular, if we take the precise representative of a least gradient function $u$, then $\partial \{ u \geq t \}$ is an analytical minimal surface for every $t$. For this reason, we will in this paper \textbf{always assume that} $u$ \textbf{is the precise representative of a least gradient function} in order to be able to state any pointwise results.

The following result is a weak maximum principle for least gradient functions, as it states that each of the level superlevel sets cannot have compact support in $\Omega$, i.e. the maximum value is attained on the boundary.

\begin{stw}\label{wmp}(\cite[Theorem 3.4]{Gor})
Let $\Omega \subset \mathbb{R}^k$, where $2 \leq k \leq 7$ and suppose $u \in BV(\Omega)$ is a function of least gradient. Then for every $t \in \mathbb{R}$ the set $\partial \{ u \geq t \}$ is empty or it is a sum of minimal surfaces $S_{t,i}$, pairwise disjoint in $\Omega$, which satisfy $\partial S_{t,i} \subset \partial \Omega$, where $\partial S_{t,i}$ is the boundary of $S_{t,i}$ in $\partial \{ u \geq t \}$.
\end{stw}

We conclude this subsection by bringing together Lemmata \ref{min+max} and \ref{traceofminimum} to notice that

\begin{wn}\label{minmax}
If $u, v \in BV(\Omega)$ are solutions to the least gradient problem with boundary data $h \in L^1(\Omega)$ in the sense of traces, then so are $\min(u,v)$ and $\max(u,v)$. \qed
\end{wn}


\subsection{Pointwise properties}

The next two Lemmata give us some insight about local form of superlevel sets of a least gradient function. Namely, locally there is only one connected component of $\{ u \geq t \}$ around any point inside $\Omega$ (this statement may obviously fail at th boundary). Secondly, absence of connected components of $\{ u \geq t \}$ passing through a given point imply that there are none in some neighbourhood of this point. 

\begin{lem}\label{lem:tylkojednowkuli}
Let $2 \leq k \leq 7$. Let $u \in BV(\Omega)$ be a least gradient function. Let $E_t = \{ u \geq t \}$ and take $x \in \partial E_t$. Then there exists a ball $B(x,r)$ such that there is only one connected component of $\partial E_t$ intersecting this ball.
\end{lem}

\begin{dd}
Suppose otherwise: for each ball $B(x,r)$ with $x \in \partial E_t \cap \Omega$ we have at least two connected components of $\partial E_t$ intersecting this ball. Let us call them $S_0$ and $S_1$. The sets $S_0 \cap \overline{B(x,r)}$ and $S_1 \cap \overline{B(x,r)}$ are compact and disjoint due to Proposition \ref{wmp}. Let $d$ be the (positive) distance between these sets. Then, by our assumption, $B(x, \frac{d}{2})$ contains at least two connected components of $\partial E_t$ and neither of them is $S_1$, thus there were at least three components intersecting $B(x,r)$; by repeating this reasoning we obtain that there are infinitely many connected components of $\partial E_t$ in each ball. 

By Proposition \ref{wmp} and the Alexander duality theorem, see \cite[Theorem 27.10]{GH}, $S_0$ divides $\Omega$ into two disjoint sets, $\Omega_+$ and $\Omega_-$, so there are infinitely many connected components in either $\Omega_+$ or $\Omega_-$; up to renumbering of $S_k$ we may assume there are infinitely many connected components of $\partial E_t$ between $S_0$ and $S_1$. Then for each $S_k$ between $S_0$ and $S_1$ the area $\mathcal{H}^{n-1}(S_k \cap B(x,r))$ is bounded from below; let $\Pi$ by a hyperplane tangent to $S_0$ at $x$. Then the orthogonal projection of $S_k \cap B(x,r)$ onto $\Pi$ contains the orthogonal projection of $S_1 \cap B(x,r)$ onto $\Pi$, as $S_k$ is between $S_0$ and $S_1$. Then

\begin{equation*}
\mathcal{H}^{n-1}(S_k \cap B(x,r)) \geq \mathcal{H}^{n-1}(pr_\Pi (S_k \cap B(x,r))) \geq \mathcal{H}^{n-1}(pr_\Pi (S_1 \cap B(x,r))) > 0, 
\end{equation*}
so the total variation of $D \chi_{E_t}$ is infinite: 

\begin{equation*}
|D \chi_{E_t}|(\Omega) \geq |D \chi_{E_t}|(B(x,r)) = \sum_{k = 0}^{\infty} \mathcal{H}^{n-1}(S_k \cap B(x,r))  = +\infty.
\end{equation*}
As by Theorem \ref{twierdzeniezbgg} $\chi_{E_t}$ is a function of least gradient, we have reached a contradiction. \qed
\end{dd}

\begin{lem}\label{lem:zawieraniewkuli}
Let $2 \leq k \leq 7$. Suppose that $u \in BV(\Omega)$ is a function of least gradient. Let $E_t = \{ u \geq t \}$. Suppose that $x \in \Omega$ is a point of continuity of $u$, $u(x) = t$ and $x \notin \partial E_t$. Then there exists a ball $B(x,r) \subset E_t$.
\end{lem}

\begin{dd}
We have two possibilities: either $|D \chi_{E_t}|(B(x,r)) > 0$ for all $r > 0$ or for sufficiently small $r$ we have $|D \chi_{E_t}|(B(x,r)) = 0$. 

In the first case we set $d = \text{dist}(x, \partial\Omega)$ and take any $r < d$. As $|D \chi_{E_t}|(B(x,r)) > 0$ for all $r > 0$, we have at least one (and thus infinitely many) connected component of $\partial E_t$ intersecting $B(x,r)$. Now the proof follows the same lines as the proof of the previous lemma.

In the second case we take such $r$. By relative isoperimetric inequality we have either $B(x,r) \subset E_t$ or $B(x,r) \cap E_t = \emptyset$ (remember that we consider the precise representative of $u$). But the second condition cannot hold, as $u(x) = t$ and $x$ is a point of continuity.  \qed
\end{dd}

\begin{uw}
By \cite[Proposition 3.6]{Gor} it suffices to assume that $x \notin \partial E_t$ for any $t \in \mathbb{R}$; then $x$ is a point of continuity of $u$. Moreover, suppose that $u \in BV(\Omega)$ is a function of least gradient, $u(x) = t$, $x \notin \partial \{ u \geq t \}$ for any $t$ and $x \notin \partial \{ u \leq t \}$ for any $t$. Then there exists a ball $B(x,r) \subset \{ u = t \}$.
\end{uw}

\begin{prz}\label{prz:nieciaglynaotoczeniu}
However, if $x \notin \partial E_t$ for any $t \in \mathbb{R}$, this does not mean that $u$ is continuous in any open neighbourhood of $x$; take a nonincreasing function on $[-1,1]$ (one can easily produce identical examples on $\mathbb{R}^N$) defined by the formula

\begin{equation*}
u(x) = \twopartdef{2^{\lfloor \frac{1}{x} \rfloor}}{x < 0}{0}{x \geq 0}.
\end{equation*}
This function is continuous at $0$, yet it is not continuous on any open interval $(-\delta, \delta)$. We see that $0 \notin \partial \{ u \geq t \}$ for any $t$; for $t \leq 0$ it is impossible, as on the whole domain the function is nonnegative. For $t > 0$ we will find $x_0 \in (-1,0)$ such that for $x > x_0$ we have $u(x) < t$. Thus, as the function is not constant anywhere near $0$, by the previous remark $0 \in \partial \{ u \leq 0 \}$.
\end{prz}

The following result states that there can be only countably many $t$ such that $\partial \{ u > t \} \neq \partial \{ u \geq t \}$.

\begin{lem}\label{lem:duzepoziomice}
Let $2 \leq k \leq 7$. Suppose that $u \in BV(\Omega)$ is a function of least gradient. We have $\partial \{ u > t \} \neq \partial \{ u \geq t \}$ if and only if $|\{u = t \}| > 0$.
\end{lem}

\begin{dd}
Suppose that $\partial \{ u > t \} \neq \partial \{ u \geq t \}$. Obviously $\{ u > t \} \subset \{ u \geq t \}$ and by Theorem \ref{twierdzeniezbgg} their boundaries are minimal surfaces. We have two possibilities: either there is a connected component $S$ of $\partial \{ u \geq t \}$ such that $S \cap \partial \{ u > t \} = \emptyset$ or there is not. In the first case we easily see, for example using Lemma \ref{lem:tylkojednowkuli}, that $|\{ u \geq t \} \backslash \{ u > t \}| = |\{ u = t \}| > 0$. In the second case, let us see that by \cite[Theorem 2.2]{SWZ}, later stated as Proposition \ref{stw:sternbergpowmin}, if a connected component of $\{ u \geq t \}$ and a connected component of $\{ u > t \}$ intersect, then they are equal; thus the second case cannot happen. 

we have either $\partial \{ u > t \} = \partial \{ u \geq t \}$ or for some connected component $S$ of $\partial \{ u \geq t \}$ we have $S \cap \partial \{ u > t \} = \emptyset$. 

In the other direction, suppose that $|\{ u = t \}| > 0$. Take a point $x \in \partial \{ u = t \} \cap \Omega \subset (\partial \{ u \geq t \} \cup \{ u \leq t \}$ (if we omitted the intersection with $\Omega$, we would have an additional summand $\partial\Omega$ on the RHS of the inclusion). By the previous remarks, as $u$ is not constant in any neighbourhood of $x$, we have exactly one of the sets $\partial \{ u \geq t \}$ and $\{ u \leq t \}$ passing through $x$.
\end{dd}


\subsection{The weak maximum principle}

Unfortunately, the weak maximum principle as presented in Proposition \ref{wmp} is not enough for our considerations. The next two results are improvements of the weak maximum principle which consider the geometry of the superlevel sets of a least gradient function near the boundary of $\Omega$: they state that two connected components of a superlevel set cannot intersect even on $\partial\Omega$. The first result is two-dimensional and it serves as a toy model for the second one, which covers the general case. Note that these results require additionally convexity of $\Omega$; however, they do not require strict convexity of $\Omega$.

\begin{lem}\label{lem:slabazm2d}
Let $\Omega \subset \mathbb{R}^2$ be a convex set with Lipschitz boundary and suppose $u \in BV(\Omega)$ is a function of least gradient. Let $E_t = \{ u \geq t \}$. Then for every $t \in \mathbb{R}$ for every point $x \in \partial\Omega$ there is at most one interval belonging to $\partial E_t$ which ends at $x$.
\end{lem}

\begin{dd}
Suppose we have at least two intervals in $\partial E_t$: $\overline{xy}$ and $\overline{xz}$. We have two possibilities: there are countably many intervals in $\partial E_t$, which end in $x$, with the other end lying in the arc $\overline{yz} \subset \partial\Omega$ which does not contain $x$; or there are finitely many. In the first case, as $\Omega$ is convex (but not necessarily stricly convex), we see that $|D \chi_{E_t}|(\Omega) = + \infty$: each of these intervals projects orthogonally onto the altitude of the triangle $xyz$ passing through $x$, so their lengths are bounded from below. Thus $\chi_{E_t} \notin BV(\Omega)$; but this contradicts Theorem \ref{twierdzeniezbgg}.

Now we move to the first case. If there are finitely many such intervals, then without loss of generality we may assume that $\overline{xy}$ and $\overline{xz}$ are adjacent. This situation is depicted on Figure 1 on the left hand side. Consider the function $\chi_{E_t}$. In the area enclosed by the intervals $\overline{xy}, \overline{xz}$ and the arc $\overline{yz} \subset \partial\Omega$ not containing $x$ we have $\chi_{E_t} = 1$ and $\chi_{E_t} = 0$ on the two sides of the triangle (or the opposite situation, which we handle similarly). Then $\chi_{E_t}$ is not a function of least gradient: the function $\widetilde{\chi_{E_t}} = \chi_{E_t} - \chi_{\Delta xyz}$ has strictly smaller total variation due to the triangle inequality. This again contradicts Theorem \ref{twierdzeniezbgg}. \qed
\end{dd}

\begin{figure}[h]\label{fig:weakmp}
\caption{Weak maximum principle}
\includegraphics[scale = 0.20]{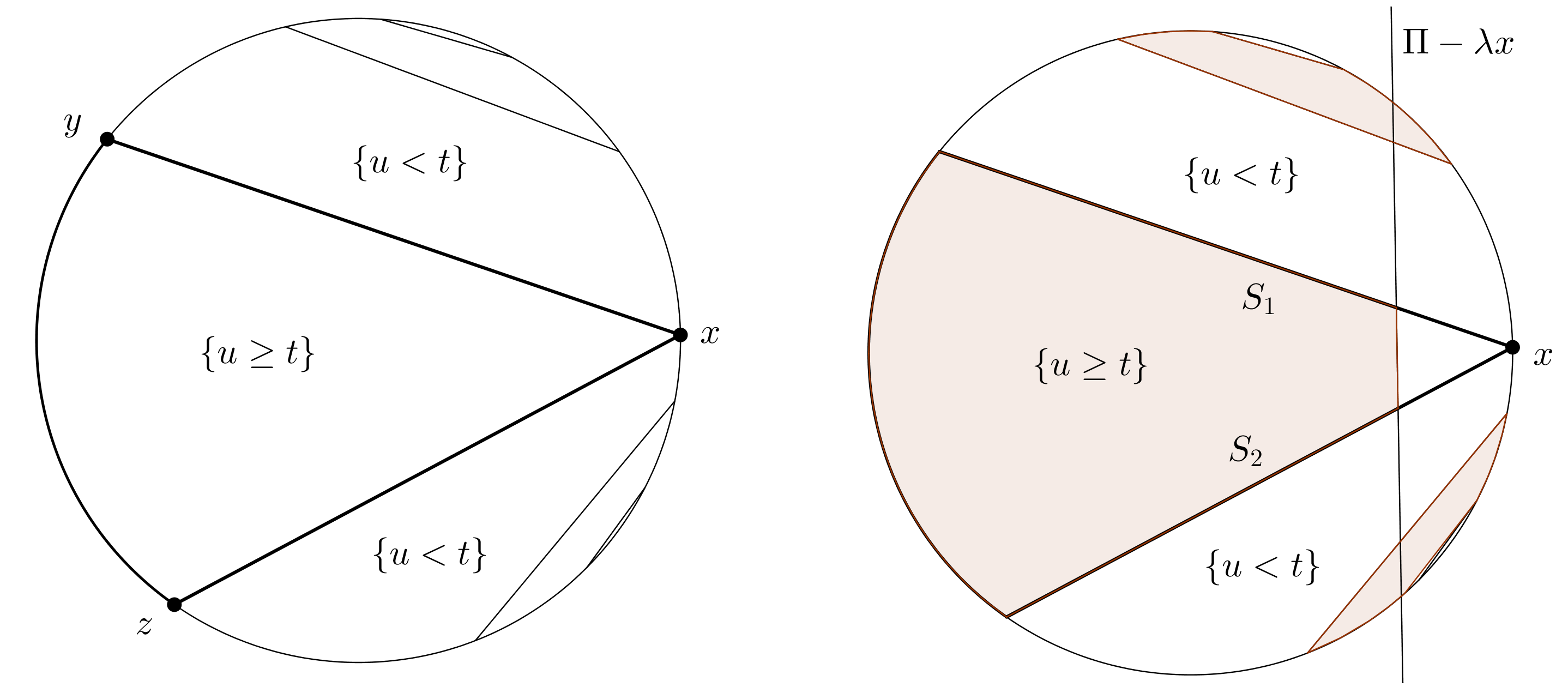}
\end{figure}

In the more general case, we have to state the result and its proof more carefully. There are two main reasons: firstly, an interval divides $\Omega$ into two simply-connected open sets, what may fail in higher dimensions: for a simple example, consider $\Omega$ to be a ball in $\mathbb{R}^3$ and $\partial \{ u \geq t \}$ to be a catenoid. Secondly, we may not use the triangle inequality and we have to rely on projections, so the geometrical part becomes more complicated. 

\begin{stw}\label{stw:slabazasadamaksimumplus}
Let $\Omega \subset \mathbb{R}^k$, where $2 \leq k \leq 7$, be a convex set with Lipschitz boundary and suppose $u \in BV(\Omega)$ is a function of least gradient. Then for every $t \in \mathbb{R}$ the boundary of the set $E_t = \{ u \geq t \}$ is a sum of minimal surfaces $S_{t,i}$, without self-intersections, with closures pairwise disjoint in $\overline{\Omega}$. 
\end{stw}

\begin{dd}
We only have to prove that intersection does not take place on $\partial \Omega$. Let $x \in \overline{S_1} \cap \overline{S_2} \cap \partial \Omega$, where $S_1$ and $S_2$ are two different connected components of $\partial E_t$. We know that $S_1$ divides $\Omega$ into two disjoint (but not necessarily connected) parts, $\Omega_1^+$ and $\Omega_1^-$; similarly $S_2$ divides $\Omega$ into $\Omega_2^+$ and $\Omega_2^-$. Among these, due to Proposition \ref{wmp}, there is only one set of the form $\Omega_1^\pm \cap \Omega_2^\pm$, which lies between $S_1$ and $S_2$, i.e. has both of these sets as parts of its boundary. Without loss of generality it is $\Omega_1^+ \cap \Omega_2^+$. 

If there is another minimal surface $S_3 \subset \partial E_t$ such that $x \in \overline{S_3}$ and $S_3 \subset \Omega_1^+ \cap \Omega_2^+$, then we may replace $S_2$ by $S_3$; this way we can assume that $S_1$ and $S_2$ are adjacent minimal surfaces (the case that there are countably many connected components is excluded similarly as in the proof of Lemma \ref{lem:tylkojednowkuli}). Without loss of generality $\Omega_1^+ \cap \Omega_2^+ \subset E_t$ and it is a connected component of $E_t$.

Consider the hyperplane $\Pi$ tangent to $\partial \Omega$ at $x$ (as $\partial\Omega$ is only Lipschitz, such a hyperplane might not exist; in that case take any of the supporting hyperplanes). Theorem \ref{twierdzeniezbgg} implies that $\chi_{E_t}$ is a function of least gradient in $\Omega$. Now consider a competitor $\chi_F$ constructed in the following way:

- in $\Omega \backslash (\Omega_1^+ \cap \Omega_2^+)$ we have $F = E_t$;

- there are two subsets of $\Omega$ bounded by $\partial\Omega$ and $\Pi$ translated by $\lambda x$ for sufficiently small $\lambda$ (chosen with respect to $S_1, S_2$). Let $G$ be the one such that $x$;

- in $\Omega_1^+ \cap \Omega_2^+$ we take $F = E_t \backslash G$.

This situation is presented on Figure 1 on the right hand side. Here, the set $F$ is the shaded region. The characteristic function $\chi_F$ constructed this way obviously satisfies $T\chi_{E_t} = T\chi_F$. Moreover, let us see that

$$|D\chi_F|(\Omega) = |D\chi_F|(\Omega \backslash (\Omega_1^+ \cap \Omega_2^+)) + \mathcal{H}^{k-1}(S_1 \cap (\Omega \backslash G)) +  $$
$$ + \mathcal{H}^{k-1}(S_2 \cap (\Omega \backslash G)) + \mathcal{H}^{k-1}((\Pi - \lambda x) \cap (\Omega_1^+ \cap \Omega_2^+)) < $$
$$ < |D\chi_{E_t}|(\Omega \backslash (\Omega_1^+ \cap \Omega_2^+)) + \mathcal{H}^{k-1}(S_1) + \mathcal{H}^{k-1}(S_2) = |D\chi_{E_t}|(\Omega), $$
as the first summands are the same and projection onto $(\Pi - \lambda x)$ delivers strict inequality in the remaining summands. We have reached a contradiction with Theorem \ref{twierdzeniezbgg}. \qed
\end{dd}

Finally, let us see that convexity of $\Omega$ in Lemma \ref{lem:slabazm2d} and Proposition \ref{stw:slabazasadamaksimumplus} cannot be relaxed.

\begin{prz}
Denote by $\varphi$ the angular coordinate in the polar coordinates on the plane. Let $\Omega = B(0,1) \backslash (\{ \frac{\pi}{4} \leq \varphi \leq \frac{3\pi}{4} \} \cup \{ 0 \}) \subset \mathbb{R}^2$, i.e. the unit ball with one quarter removed. Note that the set $\Omega$ is star-shaped, but it is not convex. Take the boundary data $f \in L^1(\partial \Omega)$ to be

$$ f(x,y) = \twopartdef{1}{y \geq 0}{0}{y < 0.}$$
Then the solution to the least gradient problem is the function (defined inside $\Omega$)

$$ u(x,y) = \twopartdef{1}{y \geq 0}{0}{y < 0,}$$
in particular $\partial \{ u \geq 1 \}$ consists of two horizontal intervals whose closures intersect on $\partial \Omega$ at the point $(0,0)$. \qed
\end{prz}


\section{Uniqueness}

This section is devoted to proving the main result of this paper, namely uniqueness of solutions of the least gradient problem except for a set where the solution is locally constant. The proof is valid in dimensions up to seven, i.e. such that boundaries of superlevel sets are analytic minimal surfaces. However, much of the proof is simplified in the planar case, i.e. when $k = 2$. In the beginning, let us underline the fact that {\bf we are always dealing with exact representatives} of least gradient functions, and thus we may discuss pointwise properties of least gradient functions. Our main tools will be Theorem \ref{twierdzeniezbgg}, connecting least gradient functions to minimal surfaces, and the following variant of the maximum principle for minimal graphs:

\begin{stw}\label{stw:sternbergpowmin} $($\cite[Theorem 2.2]{SWZ}) \\
Suppose $E_1 \subset E_2$ and let $\partial E_1, \partial E_2$ are area-minimizing in an open set $U$. Further, suppose $x \in \partial E_1 \cap \partial E_2 \cap U$. Then $\partial E_1$ and $\partial E_2$ coincide in a neighbourhood of $x$. \qed
\end{stw}

Let us note that $\partial E_1$ and $\partial E_2$ agree on their respective connected components. Now we recall the statement of Theorem \ref{uniqueness}: \vspace{2mm} \\
\textbf{Theorem \ref{uniqueness}.} \textit{
Let $\Omega \subset \mathbb{R}^k$, where $2 \leq k \leq 7$, be an open bounded convex set with Lipschitz boundary. Let $u, v$ be precise representatives of functions of least gradient in $\Omega$ such that $Tu = Tv = h$. Then $u = v$ on $\Omega \backslash (C \cup N)$, where both $u$ and $v$ are locally constant on $C$ and $N$ has Hausdorff dimension at most $k-1$.} \vspace{2mm}

For the whole section we introduce the following notation: let $u, v \in BV(\Omega)$ be two functions of least gradient with the same trace. Let $E_t = \{ u \geq t \}$ and $F_s = \{ v \geq s \}$. The proof will consist of four major steps: 

\textbf{1}. We prove that if $\partial E_t \cap \partial F_t \neq \emptyset$, then they coincide on their respective connected components; this gives a partition of $\Omega$.

\textbf{2}. We look at the structure of the set $E_t \backslash F_t$.

\textbf{3}. We use this knowledge to prove that for $t \neq s$ we have $\partial E_t \cap \partial F_s \subset J_u \cup J_v$.

\textbf{4}. We introduce a singular set $N$ with Hausdorff dimension $k-1$. We infer local properties of $u$ and $v$ from the steps above; case-by-case analysis proves uniqueness outside of $C \cup N$.

The proof is much easier to visualize in the two-dimensional case. This is most striking in Step 3 of the proof, therefore Step 3 will be proved in two stages: firstly in a two-dimensional setting with far fewer technical difficulties, secondly in the general setting with the two-dimensional proof serving as an illustration.


\textbf{Proof of Theorem \ref{uniqueness}.}


\textbf{Step 1.} \textbf{Let $x \in \partial E_t \cap \partial F_t$. Then the respective connected components of $\partial E_t$ and $\partial F_t$ coincide.}

We begin with noting that Step 1 remains the same for $2 \leq k \leq 7$. From Proposition \ref{wmp} we know that $\partial E_t$ is an at most countable sum of minimal surfaces which do not intersect inside $\Omega$ (including self-intersections). By Lemma \ref{lem:slabazm2d} (in the two-dimensional case) or by Proposition \ref{stw:slabazasadamaksimumplus} (in the general case) they do not intersect on $\partial\Omega$.

Let $w = \min(u,v)$. We assumed $Tu = Tv = h \in L^1(\partial\Omega)$. By Corollary \ref{minmax} $w$ is another function of least gradient with boundary data $h$. Consider $H_t = \{ \min(u,v) \geq t \} = E_t \cup F_t$. Let $x \in \partial E_t \cap \partial F_t$ and let $S_u$, $S_v$ be connected components of $\partial E_t$ and $\partial F_t$ respectively containing $x$.

Using Lemma \ref{lem:tylkojednowkuli} we can find a ball $B(x,r) \subset \Omega$ that intersects only $S_u$ and $S_v$ among all connected components of $\partial E_t$ and $\partial F_t$. Now we have two possibilities:

(1) For every sequence $\rho_n \rightarrow 0$ we have $S_u \cap B(x, \rho_n) \neq S_v \cap B(x, \rho_n)$. In this case every neighbourhood of $x$ intersects $\{ u, v < t \}$, thus $x \in \partial H_t$.

(2) There is an open ball $B(x, \rho)$ with $\rho < r$ such that $S_u \cap B(x, \rho) = S_v \cap B(x, \rho)$. 

Now, if there is another point $y \in \partial S_u \cap \partial S_v$ such that condition (1) holds with $x'$ in place of $x$, then $x'  \in \partial H_t$ and we may proceed to the next paragraph. If condition (2) holds for every $y \in S_u \cap S_v$, then as the intersection of two minimal surfaces is a closed set in $\Omega$, we have $S_u = S_v$.

By the reasoning above, we have $x \in \partial H_t$ (or $x' \in \partial H_t$). As $E_t \subset H_t$, by Proposition \ref{stw:sternbergpowmin} we have that $S_u = S_w$, where $S_w$ is the connected component of $\partial H_t$ containing $x$. Similarly, as $F_t \subset H_t$, we have $S_v = S_w$; thus $S_u = S_v$.


\textbf{Step 2.} \textbf{The structure of $E_t \backslash F_t$ for all but countably many $t \in \mathbb{R}$.}

By the Alexander duality theorem, see \cite[Theorem 27.10]{GH}, each of the surfaces $S_u \subset \partial E_t$ divides $\Omega$ into two open (for $k > 2$ not necessarily connected) sets $\Omega_+$ and $\Omega_-$ (in two dimensions one may use the Jordan curve theorem). If any other connected component of $E_t$ or $F_t$ intersects $\Omega_\pm$, then by Step 1 it entirely lies in $\Omega_\pm$. Now take any connected component of $\partial E_t$ or $\partial F_t$ which lies in $\Omega_\pm$ (if such exists) and it divides $\Omega_\pm$ again into two sets. This way we obtain a decomposition of $\Omega$ into at most countably many pairwise disjoint open sets $\Omega_i$; possibly dividing them into their connected components, we may assume them to be connected.

Notice that the set $E_t \backslash F_t$ may not touch the boundary of $\partial \Omega$ on a set of positive Hausdorff measure for all but countably many $t$. Indeed, if it has nonzero measure, then by the positivity of the trace functional we have $h = Tu \geq t \geq Tv = h$. Thus $h = t$ on a set of positive measure, which may happen only for countably many $t \in \mathbb{R}$. From now on, assume that $t$ is such that the level set $\{ h = t \} \subset \partial\Omega$ has zero Hausdorff measure.

Under this assumption, the boundary of $C_t$, a connected component of $E_t \backslash F_t$, cannot consist of parts of $\partial\Omega$ of positive area. Thus $\partial C_t$ is an at most countable sum of minimal surfaces $S_i \subset \partial E_t$ and $T_j \subset \partial F_t$. Let $S_i$ and $T_j$ be the connected components of $\partial C_t$ which belong to $\partial E_t$ and $\partial F_t$ respectively. We may say that $S_i$ and $T_j$ {\bf interlace}, as $\overline{S_i} \cap \overline{S_j} = \emptyset$ for $i \neq j$ due to Lemma \ref{lem:slabazm2d}, while $\overline{S_i}$ and $\overline{T_j}$ must overlap for some $i$ and $j$. One way to imagine this is, in the two-dimensional setting, that $C_t$ is a $2n-$sided polygon such that the even sides belong to $\partial E_t$ and odd sides belong to $\partial F_t$; a three-dimensional example could be $S_1$ to be a part of a vertical catenoid and $T_1$ and $T_2$ be two horizontal disks. Here, $C_t$ is the set bounded by these three surfaces.

As $\chi_{E_t}$ is a function of least gradient in $\Omega$, taking as a competitor
the function $\chi_{E_t} - \chi_{C_t}$ (note that $T\chi_{C_t} = 0$) we obtain that 

$$ \sum_{i = 1}^\infty \mathcal{H}^{n-1}(S_i) \leq \sum_{j = 1}^\infty \mathcal{H}^{n-1}(T_j).$$
Similarly, as $\chi_{F_t}$ is a function of least gradient, taking the function $\chi_{F_t} + \chi_{C_t}$ as a competitor we have

$$\sum_{j = 1}^\infty \mathcal{H}^{n-1}(T_j) \leq \sum_{i = 1}^\infty \mathcal{H}^{n-1}(S_i).$$
This implies that for every $t$ the set $C_t$ satisfies what we will call the {\bf Green's formula}, i.e. we have

\begin{equation}\label{eq:wzorgreenadlact}
\sum_{i = 1}^\infty \mathcal{H}^{n-1}(S_i) = \sum_{j = 1}^\infty \mathcal{H}^{n-1}(T_j).
\end{equation}


\textbf{Step 3: the two-dimensional case.} \textbf{For all but countably many $t, s$ such that $t \neq s$ we have $\partial E_t \cap \partial F_s \subset J_u \cup J_v$.}

Without loss of generality we have $s < t$. Let $t, s$ be as in Step 2, i.e. such that $E_t \backslash F_t$ does not touch $\partial \Omega$ on a set of positive measure, so the interlacing condition and Green's formula are satisfied. Suppose that $x \in \partial E_t \cap \partial F_s$ and that $u, v$ are continuous at $x$. Obviously $x \in E_t \backslash F_t$. Consider $C_t$, the connected component of $E_t \backslash F_t$ containing $x$. As $x \in \partial E_t \cap \partial F_s$ and $u, v$ are continuous at $x$, there is a point $y$ in the neighbourhood of $x$ such that $y \in E_s \backslash F_s$. Similarly, let $C_s$ be a connected component of $E_s \backslash F_s$ containing $y$. 

The proof of this Step is much more clear in dimension two and we may rely on triangle inequality in place of Proposition \ref{stw:sternbergpowmin}. The situation is represented on Figure 2. Step 2 implies that the polygon $C_s$ has at most countably many vertices $p_i$ and (due to interlacing condition) its sides $p_i p_{i+1}$ belong alternately of connected components of $E_t$ and $F_t$. Similarly, the polygon $C_t$ has at most countably many vertices $q_i$ and its sides $q_i q_{i+1}$ consist alternately of connected components of $E_s$ and $F_s$. Finally, the points $r_i$ are intersections between sides of the two polygons, such that $r_1$ and $r_2$ lie on $q_1 q_2$, $r_2$ and $r_3$ lie on $p_2 p_3$ and so on. The enumeration is chosen so that $x = r_1$. If there is a finite number $N_0$ of intervals, then we employ the notation that $p_1 = p_{N_0 + 1}$ and so on.

The structure of these sets (the intersection is a polygon with trapezoids belonging to $C_t$ and $C_s$ alternately) is as on Figure 2, because $E_t \subset E_s$ and $F_t \subset F_s$ (we encourage the reader to draw how do the sets $E_t$, $F_t$, $E_s$ and $F_s$ look like). The only thing that can be different is that some of the intervals may coincide when we have a jump, i.e. $q_3 q_4 = p_3 p_4$. For now, let us assume this is not the case and this will be discussed later.

\begin{figure}[h]\label{fig:etft}
\caption{The sets $C_t$ and $C_s$}
\includegraphics[scale = 0.40]{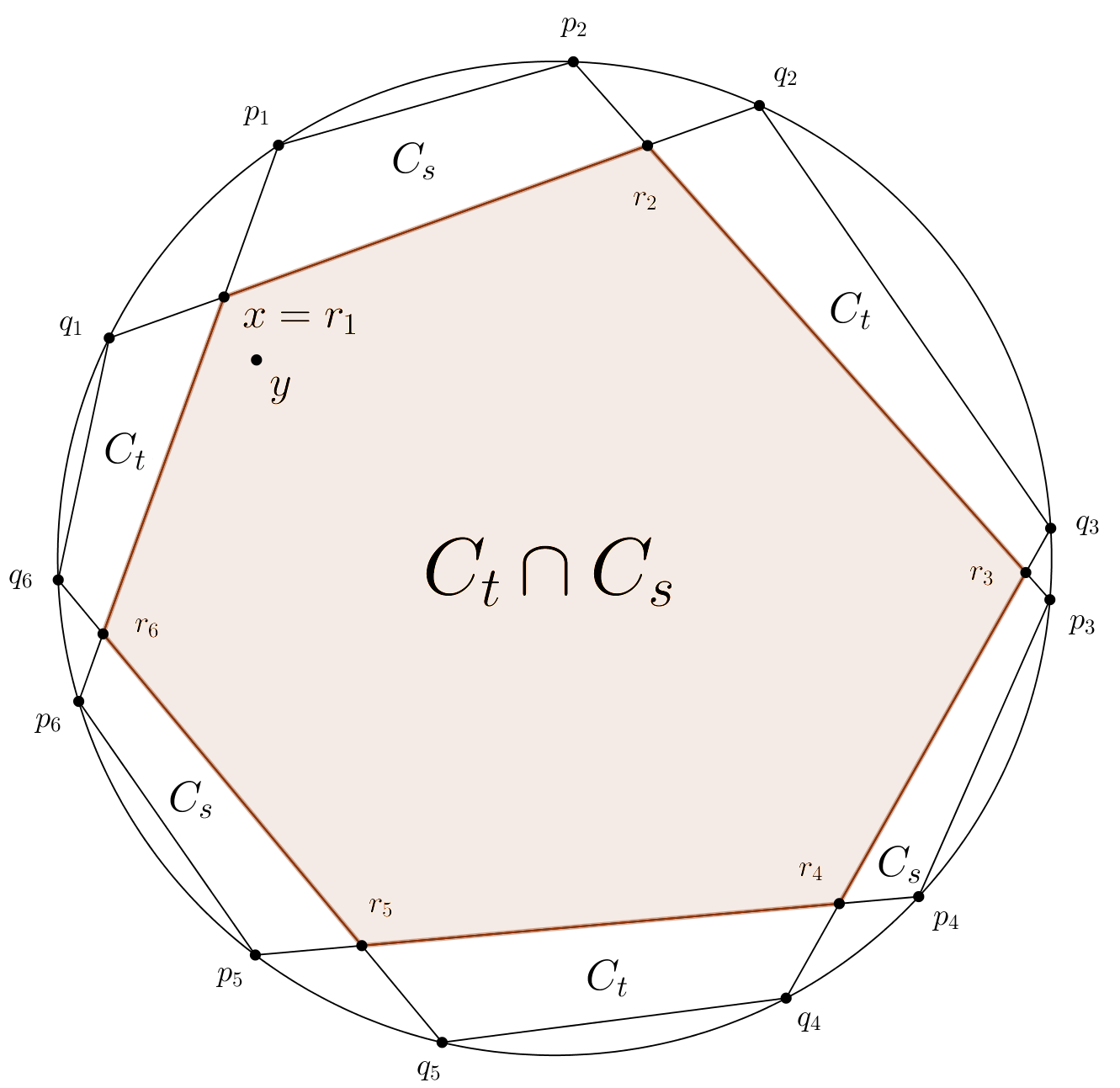}
\end{figure}

Let us look at the little trapezoids at the sides of $C_t \cap C_s$. By triangle inequality for every $i$ we have 

$$|p_{2i-1} p_{2i}| < |p_{2i-1} r_{2i-1}| + |r_{2i-1} r_{2i}| + |r_{2i} p_{2i}|;$$
$$|q_{2i} q_{2i+1}| < |q_{2i} r_{2i}| + |r_{2i} r_{2i+1}| + |r_{2i+1} q_{2i+1}|.$$
We sum up these inequalities and use the collinearity of $q_{2i-1}, r_{2i-1}, r_{2i}, q_{2i}$ and the collinearity of $p_{2i}, r_{2i}, r_{2i+1}, p_{2i+1}$ to obtain

\begin{equation}\label{eq:badineq}
\sum_i |p_{2i-1} p_{2i}| + \sum_i |q_{2i} q_{2i+1}| < \sum_i |p_{2i} p_{2i+1}| + \sum_i |q_{2i-1} q_{2i}|.
\end{equation}
This contradicts Green's formula: in the notation of Step 2, we have $S_{s,i} = p_{2i-1} p_{2i}$, $T_{s,j} = p_{2j} p_{2j+1}$, $S_{t,i} = q_{2i-1} q_{2i}$ and $T_{t,j} = q_{2j} q_{2j+1}$. Thus application of equation (\ref{eq:wzorgreenadlact}) for $C_t$ and $C_s$ implies that in equation (\ref{eq:badineq}) there should be an equality, contradiction.

Let us go back to the case where some of the intervals coincide. Then the corresponding inequality ceases to be strict. However, at least one inequality is not strict: the inequality for $i = 1$, as we assumed that there is no jump at $x$. Thus the proof still holds.


\textbf{Step 3: the general case.} \textbf{For all but countably many $t, s$ such that $t \neq s$ we have $\partial E_t \cap \partial F_s \subset J_u \cup J_v$.}

We proceed similarly to the two-dimensional case. We are going to prove the statement by contradiction: without loss of generality we have $s < t$. Again, let $t, s$ be as in Step 2, i.e. such that $E_t \backslash F_t$ does not touch $\partial \Omega$ on a set of positive measure, so the interlacing condition and Green's formula are satisfied. Suppose that $x \in \partial E_t \cap \partial F_s$ and that $u, v$ are continuous at $x$. Obviously $x \in E_t \backslash F_t$. Consider $C_t$, the connected component of $E_t \backslash F_t$ containing $x$. By Step 2 the boundary of $C_t$ consists of at most countably many minimal surfaces, $S_{t,i}$ and $T_{t,j}$, the connected components of $\partial E_t$ and $\partial F_t$ respectively. As $x \in \partial E_t \cap \partial F_s$ and $u, v$ are continuous at $x$, there is a point $y$ in the neighbourhood of $x$ such that $y \in E_s \backslash F_s$. Similarly, let $C_s$ be a connected component of $E_s \backslash F_s$ containing $x$. 

Without loss of generality assume that $x \in S_{t,1}$. This divides $\Omega$ into two parts: $\Omega_{t,1}^+$ and $\Omega_{t,1}^-$. $\Omega_{t,1}^-$ is the part of $\Omega$ which locally close to $S_{t,1}$ contains $\{ u < t \}$. If $\Omega_{t,1}^- \cap \partial E_s = \emptyset$, then $\Omega_{t,1}^- \subset E_s$; but this contradicts Step 2 for $C_s$, a connected component of $E_s \backslash F_s$. Thus there is a connected component $S_{s,1}$ of $\partial E_s$ in $\Omega_{t,1}^-$ (additionally we may pick the one closest to $S_{t,1}$). As $u$ is continuous at $x$, by Proposition \ref{stw:sternbergpowmin} $S_{s,1} \cap S_{t,1} = \emptyset$, i.e. $S_{s,1} \subset \Omega_{t,1}^-$. This reasoning mirrors the third and the last paragraph of the two-dimensional proof.

The boundary of $C_s$ contains $S_{s,1}$. Similarly to the reasoning above, using Proposition \ref{stw:sternbergpowmin} we prove that $S_{s,i} \subset \overline{\Omega_{s,i}^-}$ and $T_{t,j} \subset \overline{\Omega_{t,j}^-}$. Similarly as in the two-dimensional case, here we cannot exclude the case that $T_{t_j} = T_{s,j}$. 

Now, both $C_t$ and $C_s$ satisfy Green's formula. Explicitly, from equation (\ref{eq:wzorgreenadlact}) we have

$$ \sum_{i = 1}^\infty \mathcal{H}^{n-1}(S_{t,i}) = \sum_{j = 1}^\infty \mathcal{H}^{n-1}(T_{t,j})$$
$$ \sum_{i = 1}^\infty \mathcal{H}^{n-1}(S_{s,i}) = \sum_{j = 1}^\infty \mathcal{H}^{n-1}(T_{s,j}).$$

Let us look at $S_{s,i}$ and $T_{t,j}$, i.e. these connected components of $\partial C_t$ and $\partial C_s$ which lay outward with respect to $y$, i.e. if we draw any Jordan curve from $y$ to any point in $S_{s,i}$, then it intersects a point from $S_{t,i}$ (as illustrated in the two-dimensional case on Figure \ref{fig:etft}); similarly, if we draw any Jordan curve from $y$ to any point in $T_{t,j}$, then it intersects a point from $T_{s,j}$. Finally, we will notice that

$$ \sum_{i = 1}^\infty \mathcal{H}^{n-1}(S_{s,i}) + \sum_{j = 1}^\infty \mathcal{H}^{n-1}(T_{t,j}) < \sum_{i = 1}^\infty \mathcal{H}^{n-1}(S_{t,i}) = \sum_{j = 1}^\infty \mathcal{H}^{n-1}(T_{s,j}).$$
Take the surface $S_{s,i}$. It divides $\Omega$ into $\Omega_{s,i}^+$ and $\Omega_{s,i}^-$. As $\chi_{E_s}$ is a function of least gradient, then its localized version $\chi_{E_s \cap \Omega_{s,i}^+}$ is as well. Let $G$ be the set bounded by $S_{s,i}$, $S_{t,i}$ and these $T_{s,j}$ which intersect $S_{t_i}$. Consider a competitor $\chi_F$, where $F$ is the set $(E_s \cap \Omega_{s,i}^+) \backslash G$. The situation is presented on Figure 3, which is a zoomed-in version of Figure 2; the set $F$ is the shaded region, the set $G$ is the trapezoid on the top and $\Omega_{s,i}^+$ is everything below the interval $S_{s,i}$.

\begin{figure}[h]\label{fig:etft_zoom}
\caption{$\partial E_t \cap \partial F_s \subset J_u \cup J_v$}
\includegraphics[scale = 0.40]{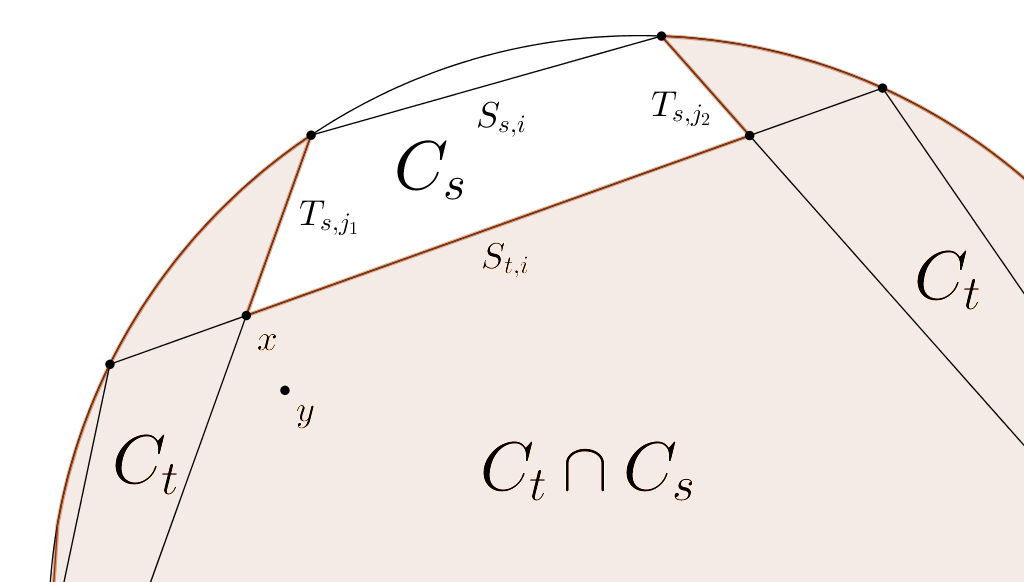}
\end{figure}

As $\chi_{E_s \cap \Omega_{s,i}^+}$ is a function of least gradient, then by comparing it to $\chi_F$ we obtain

$$ \mathcal{H}^{k-1}(S_{s,i}) \leq \sum_j \mathcal{H}^{k-1}(T_{s,j} \cap \Omega_{t,i}^-) + \mathcal{H}^{k-1}(S_{t,i} \cap \overline{G}).$$
Moreover, this inequality is strict. If it was not strict, then the surface consisting of parts of $T_{s,j}$ and $S_{t,i}$, i.e. the boundary of $G$ minus $S_{s,i}$, would be a minimal surface. But then by Proposition \ref{stw:sternbergpowmin} it equals $S_{t,i}$, as it intersects $S_{t,i}$ and $G \subset \Omega_{t,i}^-$.

This contradicts the Green's formula for $C_t$ and $C_s$, so our claim is proved. \qed


{ \bf Step 4.} Finally, we may define the set $N$. At first, recall that $J_u$ denotes the jump set of $u$ and that for any function $u \in BV(\Omega)$ we have $\dim_H J_u = k-1$. 

By Lemma \ref{lem:duzepoziomice} we have $\partial \{ u \geq t \} \neq \partial \{ u > t \} = \partial \{ u \leq t\}$ for at most countably many $t \in \mathbb{R}$. Similarly, the set $\{ h = t \} \subset \partial\Omega$ has positive Hausdorff measure for at most countably many $t$. Let us denote the (at most countable) set of $t \in \mathbb{R}$ satisfying either of these conditions by $T_u$. Let 

\begin{equation*}
B_u = \bigcup_{t \in T_u} (\partial \{ u \geq t \} \cup \partial \{ u \leq t\}).
\end{equation*}
We observe that this set has Hausdorff dimension at most $k-1$: each of the sets $\partial \{ u \geq t \}$ is a minimal surface with finite Hausdorff measure, and the set $B_u$ is an at most countable sum of such sets (as the function from Example \ref{prz:nieciaglynaotoczeniu} shows, it does not have to have finite Hausdorff measure). Now, we define a set (with Hausdorff dimension at most $k-1$)

\begin{equation*}
N = J_u \cup J_v \cup B_u \cup B_v.
\end{equation*}
Take $x \in \Omega \backslash N$. We have four possibilities: \\
1. $x \in \partial E_t \cap \partial F_t$. Then, as $u$ and $v$ are continuous at $x$, we have $u(x) = v(x) = t$. \\
2. $x \in \partial E_t \cap \partial F_s$ for $s \neq t$. This case is excluded by Step 3 of the proof. \\
3. $x \in \partial E_t$, $x \notin \partial F_s$ for any $s \in \mathbb{R}$. By Lemma \ref{lem:zawieraniewkuli} $v$ is constant on some ball around $x$ with value $s_0$. This case is excluded by the previous two points, if we consider some $s' \in (t, s_0)$. \\
4. $x \notin \partial E_t$, $x \notin \partial F_s$ for any $t, s \in \mathbb{R}$. Then by Lemma \ref{lem:zawieraniewkuli} $u,v$ are constant in some ball around $x$; thus $x \in C$.

This ends the proof of Theorem \ref{uniqueness}. \qed


\section{Classification of all solutions}

The purpose of this Section is to use Theorem \ref{uniqueness} and the knowledge obtained in Steps 1 and 2 of the proof of Theorem \ref{uniqueness} to form a complete classification of the solutions to the least gradient problem with boundary data $h \in L^1(\partial\Omega)$. We do not try to answer any questions about existence of solutions to the least gradient problem. For partial positive results, see \cite{Gor} and \cite{GRS}; for a partial negative result, see \cite{ST}. As Theorem  \ref{uniqueness} does not give us any direct information about the structure of solutions, only through comparison with another solution, we assume that at least one solution $u_0 \in BV(\Omega)$ exists and is known.

We start with a two-dimensional toy model. Then we pass to the full classification. However, the presented algorithm to find all solutions is fully applicable only in dimension two; one of the steps is to find all minimal decompositions of the set $C$, on which $u_0$ is locally constant, into sets with minimal boundary that satisfy Green's formula. This is equivalent to solving the Plateau problem, in which the spanning set is not homeomorphic to a sphere, but may fail to be connected (it may have countably many connected components) or (in dimension 4 or higher) simply-connected. Because of that, the reasoning in this section has two purposes: in dimension 2, the algorithm presented here enables us to find all the solutions; in dimensions 3 to 7, save for situations with additional symmetries, the reasoning below provides a way to determine if a function $u \in BV(\Omega)$ is a solution to the least gradient problem with boundary data $h$ without directly calculating the total variation.

\subsection{Detailed example of (non)uniqueness}

Take $\Omega = B(0,1) \subset \mathbb{R}^2$. Let $h$ be a function with six discontinuity points $p_1, ..., p_6 \in \partial \Omega$. On each of the arcs $(p_1, p_2), (p_3, p_4)$ and $(p_5, p_6)$ this function is continuous and strictly convex. It has a single minimum with value $-2$ in each of these intervals and limits equal to $-1$ at each end of these intervals. Similarly, $h$ is continuous and strictly concave on each of the intervals $(p_2, p_3), (p_4, p_5)$ and $(p_6, p_1)$. It has a single maximum with value $2$ in each of these intervals and limits equal to $1$ at each end of these intervals. It is easy to see that the function $u_0$ as on the left hand side of Figure 4 is a solution of the least gradient problem (for example by proceeding as in the proof of \cite[Theorem 4.6]{Gor}, i.e. using approximations to the boundary data and the Sternberg-Williams-Ziemer construction).

\begin{figure}[h]
\caption{Comparison of $u_0$ and $u$}
\includegraphics[scale = 0.23]{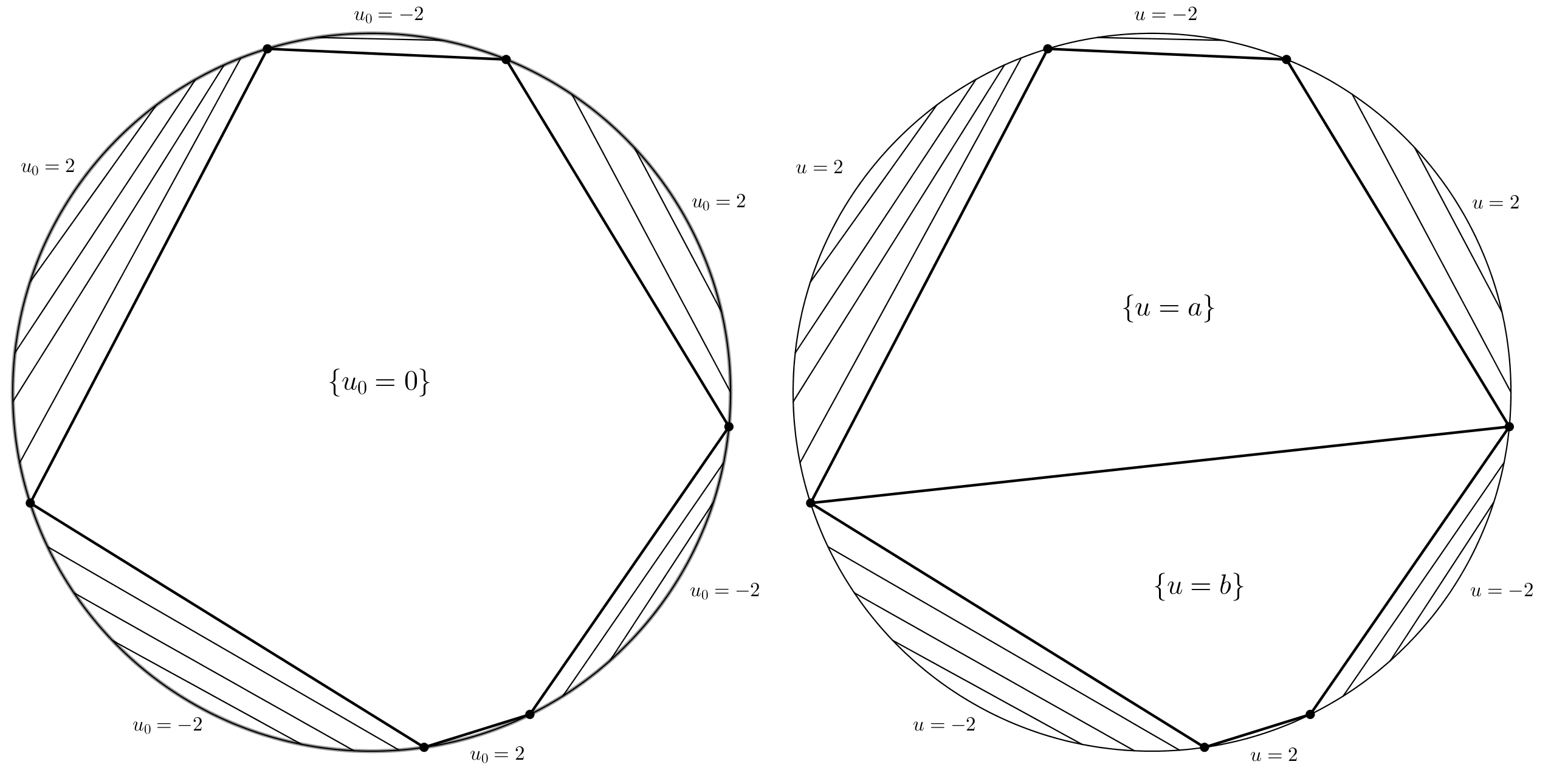}
\end{figure}

The set $C$ from the statement of Theorem \ref{uniqueness} is the hexagon $H = p_1 p_2 p_3 p_4 p_5 p_6$. Let $u \in BV(\Omega)$ be a candidate for another solution to the least gradient problem with boundary data $h$. By Theorem \ref{uniqueness} we have $u = u_0$ in $\Omega \backslash H$. We also know that $u$ is locally constant on $H$.

Let $H = \bigcup_i H_i$ such that each of the sets $H_i$ is connected and that $u$ is constant and equal $t_i$ on $H_i$. Then $\partial H_i \subset \partial \{ u \geq t_i \} \cup \partial \{ u \leq t_i \}$; by the weak maximum principle (Proposition \ref{wmp}) $\partial H_i$ composes of pairwise disjoint intervals with endpoints in $\partial\Omega$. By Proposition \ref{stw:sternbergpowmin} these intervals cannot intersect the set $\Omega \backslash H$, as they would intersect transversally some interval of the form $\partial \{ u \geq t \}$ for $t \neq t_i$. Thus these intervals have endpoints in the set $\{ p_1, ..., p_6 \}$. Moreover, analysis as in Step 2 of the proof of Theorem \ref{uniqueness} shows that the sides of the polygon $H_i$ interlace, i.e. belong alternately to $\{ u \geq t_i \}$ and $\{ u \leq t_i \}$ and satisfy Green's formula.

This means that finding all functions $u$ of least gradient with boundary data $h$ boils down to finding all subpolygons of $H$ which satisfy Green's formula. If there are none (i.e. when the hexagon is equilateral), then $u$ is constant on $H$. After a quick calculation we obtain the value $u(H)$:

$$|Du|(\Omega) = |Du|(\Omega \backslash H) + (\mathcal{H}^1(\overline{p_1 p_2}) + \mathcal{H}^1(\overline{p_3 p_4}) + \mathcal{H}^1(\overline{p_5 p_6})) |-1 - u(H)| + $$
$$ + (\mathcal{H}^1(\overline{p_2 p_3}) + \mathcal{H}^1(\overline{p_4 p_5}) + \mathcal{H}^1(\overline{p_6 p_1})) |1 - u(H)| + 0$$
and
$$|Du_0|(\Omega) = |Du_0|(\Omega \backslash H) + (\mathcal{H}^1(\overline{p_1 p_2}) + \mathcal{H}^1(\overline{p_3 p_4}) + $$
$$+ \mathcal{H}^1(\overline{p_5 p_6}) + (\mathcal{H}^1(\overline{p_2 p_3}) + \mathcal{H}^1(\overline{p_4 p_5}) + \mathcal{H}^1(\overline{p_6 p_1})) |1 - 0| + 0.$$
using Green's formula for $H$ and the fact that $u = u_0$ on $\Omega \backslash H$ we easily see that these two numbers are equal, i.e. $u$ is a function of least gradient, iff $u(H) \in [-1,1]$.
 
However, there may exist subpolygons of $H$ which satisfy Green's formula; the only possible case is two trapezoids $H_1 = p_1 p_4 p_3 p_2$ and $H_2 = p_1 p_4 p_5 p_6$ satisfying Green's formula with one common side (without loss of generality the common side is $p_1 p_4$). This situation is presented on Figure 4 on the right hand side. Let $a$ be the value on $H_1$ and $b$ the value on $H_2$. Suppose that $a \neq b$, so the situation is different from the above. A calculation similar to the one above shows that $a, b \in [-1,1]$; the only remaining problem is whether $a$ or $b$ is larger. This follows from Step 2 of the proof of Theorem \ref{uniqueness}; the sides of $H_1$ have to interlace, i.e. belong alternately to $\partial \{ u \geq a \}$ and $\{ u \leq a \}$. Thus, as $p_1 p_2 \subset \{ u \geq -1 \}$, then also $p_1 p_2 \subset \{ u \geq a \}$; but this implies that $p_1 p_4 \subset \{ u \leq a \}$, so $a < b$. Quick calculation using Green's formula for $H_1$ and $H_2$ shows that a function $u$ such that $u = u_0$ on $\Omega \backslash H$, $u = a$ on $H_1$, $u = b$ on $H_2$, $-1 \leq a \leq b \leq 1$ is of least gradient. Thus we have classified all solutions to the least gradient problem with boundary data $h$.


\subsection{Full description}

We want to find all functions of least gradient with prescribed boundary data $h \in L^1(\partial \Omega)$. Direct use of Theorem \ref{uniqueness} shows that $u = u_0$ in $\Omega \backslash C$. We want to find all admissible decompositions of $C$ into sets $C_i$ and admissible (constant) values $t_i$ of $u$ on $C_i$.

{\bf Assumption.} For simplicity, we will assume that the function $h$ has no level sets of positive measure. At the end of this chapter, we will modify this reasoning to account for such sets. Furthermore, we may assume that the set $C$, on which $u_0$ is locally constant, is connected; otherwise we could perform the same analysis for each of its connected components.

Using the reasoning from Step 2 of the proof of Theorem \ref{uniqueness} we see that $\partial C_i$ consists of an at most countable family of minimal surfaces. They belong either to $\partial \{ u \geq t_i \}$ or $\partial \{ u \leq t_i \}$ and interlace, i.e. if $S_j \subset \partial \{ u \geq t_i \}$ is a connected component of $\partial C_i$, then it intersects on the boundary with some surface $T_j \subset \partial C_i$, which is a connected component of $\partial \{ u \leq t_i \}$; furthermore, by the weak maximum principle (Proposition \ref{stw:slabazasadamaksimumplus}) it does not intersect (in $\overline{\Omega}$) any other connected component of $\partial \{ u \geq t_i \}$. The reasoning from Step 2 also implies that $C_i$ satisfies Green's formula.

Let $C = \bigcup_i C_i$ be a minimal decomposition of the set $C$ into sets $C_i$, with boundary consisting of minimal surfaces, which satisfy Green's formula and interlacing condition. We do not claim that such decomposition is unique; by minimal we only mean that no set $C_i$ can be decomposed further into multiple parts safisfying the assumptions above. Furthermore, we will denote the connected components of $\Omega \backslash C$ are by $U_i$. The trace of $u$ on $\partial U_l \cap \partial S_i$ from $U_l$ is constant (by an easy application of Lemma \ref{lem:tylkojednowkuli}) and denoted by $\alpha_{il}$.

Let us see that similarly as in the proof of \cite[Theorem 3.8]{Gor}, as $\partial C_i$ consists of minimal surfaces which provide a decomposition of $\Omega$, we may form a graph where $C_i$ are vertices and they are connected by an edge iff $\mathcal{H}^{k-1}(\partial C_i \cap \partial C_j) > 0$. This graph is a tree, i.e. it is connected (as $C$ was connected) and there is exactly one path connecting two given vertices. This time we want our graph to be directed: whenever $t_i \geq t_j$ (for neighbouring $C_i, C_j$), we draw an arrow from $C_i$ to $C_j$. In particular, if $t_i = t_j$, we draw an arrow in both directions.

The following Proposition provides a necessary and sufficient condition for given $u \in BV(\Omega)$ to be a function of least gradient with the same trace as another given function of least gradient $u_0$. 

\begin{stw}\label{stw:jakwygladaroz}
Let $\Omega \subset \mathbb{R}^k$, where $2 \leq k \leq 7$ be an open bounded convex set with Lipschitz boundary. Suppose that $h \in L^1(\partial\Omega)$ and there is at least one solution $u_0 \in BV(\Omega)$ to the least gradient problem. Then the class of solutions of least gradient problem with boundary data $h$ contains precisely the functions $u$ such that $u = u_0$ in $\Omega \backslash C$ and that $u$ has constant value $t_i$ on $C_i$ such that the following conditions are satisfied: \\
(1) In the notation introduced above, the graph for $u$ is the following: the arrows from leaves (1st level) to their neighbours (2nd level) are well defined using the interlacing condition and the same as for $u_0$. They, using the same technique, define arrows on all other edges. Then we may possibly add some arrows in the other directions (i.e. equalities $t_i = t_j$). Such graphs are possible, as there exists a graph for $u_0$; \\
(2) Whenever $\mathcal{H}^{k-1}(\partial C_i \cap \partial U_l) > 0$ and $t_i^0 \geq \alpha_{il}$, then $t_i \geq \alpha_{il}$; \\
(3) Whenever $\mathcal{H}^{k-1}(\partial C_i \cap \partial U_l) > 0$ and $t_i^0 \leq \alpha_{il}$, then $t_i \leq \alpha_{il}$.
\end{stw}

\begin{dd}
Fix any decomposition $C_i$ of the set $C$ into sets with boundary consisting of minimal surfaces, which satisfy the interlacing condition and Green's formula. Different decompositions will give us different functions of least gradient. By Theorem \ref{uniqueness} every other solution $u$ satisfies $u = u_0$ on $\Omega \backslash C$. As $u_0$ is of least gradient, $u$ is of least gradient iff $|Du|(\Omega) = |Du_0|(\Omega)$. We calculate $|Du|(\Omega)$:

\begin{equation*}
|Du|(\Omega) = |Du|(\Omega \backslash C) + \sum_{i,l} |\alpha_{il} - t_i| \mathcal{H}^{k-1}(\partial C_i \cap \partial U_l) + \sum_{i > j} |t_i - t_j| \mathcal{H}^{k-1}(\partial C_i \cap \partial C_j).
\end{equation*}
We may write here a sum over all $i, j, l$ because if $C_i$ and $C_j$ or $U_l$ do not share a boundary, the corresponding value is zero. We obtain an analogous result for $u_0$.

{\bf Sufficiency of conditions $(1) - (3)$.}  To summarize, $u$ satisfies the same inequalities between values of $u$ on $C_i$ and $C_j$ as $u_0$ and that whenever $\mathcal{H}^{k-1}(\partial C_i \cap U_l) > 0$ we have inequalities of the form

\begin{equation*}
\min_{l: \partial C_i \cap \partial U_l \subset \partial \{ u_0 \geq t_i^0 \}} \alpha_{il} \leq t_i^0 \leq \min_{l: \partial C_i \cap \partial U_l \subset \partial \{ u_0 \leq t_i^0 \}} \alpha_{il}.
\end{equation*}
We shall see that every $u$ which satisfies these properties is a function of least gradient. Denote by $s(C_i, \Omega_k)$ the function encoding inequalities between $t_i$ and $\alpha_{ik}$: let $s(C_i, \Omega_k) = 1$ if $t_i \geq \alpha_{ik}$ and $s(C_i, \Omega_k) = 0$ if the opposite inequality holds. Similarly we define $s(C_i, C_j)$. To prove that $u$ is of least gradient we have to check that $|Du|(\Omega) - |Du_0|(\Omega) = 0$.

\begin{equation*}
|Du|(\Omega) - |Du_0|(\Omega) = \sum_{i,l} |\alpha_{il} - t_i| \mathcal{H}^{k-1}(\partial C_i \cap \partial U_l) - \sum_{i,l} |\alpha_{il} - t_i^0| \mathcal{H}^{k-1}(\partial C_i \cap \partial U_l) +
\end{equation*}
\begin{equation*}
+ \sum_{i > j} |t_i - t_j| \mathcal{H}^{k-1}(\partial C_i \cap \partial C_j) - \sum_{i > j} |t_i^0 - t_j^0| \mathcal{H}^{k-1}(\partial C_i \cap \partial C_j) =
\end{equation*}
\begin{equation*}
= \sum_{i,l} (-1)^{s(C_i,U_l)}(\alpha_{il} - t_i) \mathcal{H}^{k-1}(\partial C_i \cap \partial U_l) - \sum_{i,l} (-1)^{s(C_i,U_l)}(\alpha_{il} - t_i^0) \mathcal{H}^{k-1}(\partial C_i \cap \partial U_l) +
\end{equation*}
\begin{equation*}
+ \sum_{i > j} (-1)^{s(C_i,C_j)}(t_j - t_i) \mathcal{H}^{k-1}(\partial C_i \cap \partial C_j) - \sum_{i > j} (-1)^{s(C_i, C_j)} (t_j^0 - t_i^0) \mathcal{H}^{k-1}(\partial C_i \cap \partial C_j) =
\end{equation*}
\begin{equation*}
= \sum_{i,l} (-1)^{s(C_i,U_l)}(t_i^0 - t_i) \mathcal{H}^{k-1}(\partial C_i \cap \partial U_l) + \sum_{i > j} (-1)^{s(C_i,C_j)}(t_i^0 - t_i - t_j^0 + t_j) \mathcal{H}^{k-1}(\partial C_i \cap \partial C_j) =
\end{equation*}
\begin{equation*}
= \sum_{i} (\sum_l (-1)^{s(C_i,U_l)}(t_i^0 - t_i) \mathcal{H}^{k-1}(\partial C_i \cap \partial U_l) + \sum_{j} (-1)^{s(C_i,C_j)}(t_i^0 - t_i) \mathcal{H}^{k-1}(\partial C_i \cap \partial C_j)) =
\end{equation*}
\begin{equation*}
= \sum_{i} (t_i^0 - t_i) (\sum_l (-1)^{s(C_i,U_l)} \mathcal{H}^{k-1}(\partial C_i \cap \partial U_l) + \sum_{j} (-1)^{s(C_i,C_j)} \mathcal{H}^{k-1}(\partial C_i \cap \partial C_j)) = 0,
\end{equation*}
because for every $i$ the last summand is precisely Green's formula for the sides of $C_i$. Thus every $u$ satisfying the assumptions above is a function of least gradient.

{\bf Necessity of conditions (1)-(3).} Let $C_i$ be a leaf, i.e. $C_i$ shares a boundary with only one $C_j$. Then $C_i$ shares a boundary with at least three sets of the form $U_l$. On the set $C_i$ the function $u_0$ has constant value $t_i^0$ and $u$ has constant value $t_i$. Without loss of generality assume that $\partial C_i \cap U_{l_1} \subset \partial \{ u_0 \geq t_i^0 \}$; in particular, we have $t_i^0 \geq \alpha_{i l_1}$. Using the interlacing condition we have that $\partial C_i \cap U_{l_2} \subset \partial \{ u_0 \leq t_i^0 \}$; thus $t_i^0 \leq \alpha_{i l_2}$.

Suppose that the structure of $u$ is different than the structure of $u_0$, i.e. $\partial C_i \cap U_{l_1} \subset \partial \{ u_0 \leq t_i^0 \}$. In particular $t_i \neq t_i^0$. Repeating the reasoning above we obtain that $t_i \leq \alpha_{i l_1}$ and $t_i \geq \alpha_{i l_2}$. Putting these results together, we obtain

$$t_i \leq \alpha_{i l_1} \leq t_i^0 \leq \alpha_{i l_2} \leq t_i.$$
Thus $t_i = t_i^0$, contradiction. Thus on every leaf conditions (1)-(3) are necessary. Once we do this for all the leaves, we eliminate all the leaves from the graph and repeat, treating the leaves the same as the sets $U_l$. Thus conditions $(1)-(3)$ are necessary for every $C_i$.
\qed
\end{dd}

{\bf Relaxing the assumption.} Suppose that $h$ is constant and equal to $t$ on a set $\Gamma$ of positive measure on $\partial\Omega$. Denote by $\Omega_\Gamma$ the flap enclosed by $\Gamma$; in dimension two this is particularly easy, as when $\Gamma$ is an arc, then it is a flap enclosed by $\Gamma$ and the interval connecting its endpoints. In the general case we have to remember that $\partial \{ u > t \} = \partial \{ u \leq t \}$ and $\partial \{ u \geq t \}$ compose of minimal surfaces; thus $\Gamma$ spans a set composing of minimal surfaces. Denote by $\Omega_\Gamma$ the set enclosed by these surfaces and $\partial\Omega$.

Then, by Lemma \ref{lem:zawieraniewkuli} we observe that on $\Omega_\Gamma$ the value of $u$ has to be constant and equal $t$. This value is fixed, so from now on we may treat $\Omega_\Gamma$ as one of the sets $U_l$ in the reasoning above. Thus we do not need to assume that $h$ does not have level sets of positive measure.

Let us note that Proposition \ref{stw:jakwygladaroz} has algorythmic value in case when $\Omega \subset \mathbb{R}^2$, as the only minimal surfaces are intervals, and when the decomposition into $C_i$ is finite. Finally, the following well-known examples serve as an illustration to this result:

\begin{prz}
Let $\Omega = B(0,1) \subset \mathbb{R}^2$. \\
(1) $h$ has a single maximum and a single minimum and $\partial \Omega$ can be divided into two arcs, on which $h$ is monotone. Then the solution to the least gradient problem is unique; \\
(2) $h$ takes only three values: $0$ on the arc $(p_1, p_2)$, $\alpha_1 > 0$ on the arc $(p_2, p_3)$, and $\alpha_1 + \alpha_2 > \alpha_1$ on the arc $(p_3, p_1)$ (see \cite[Section 3.4]{GRS}). Then the solution to the least gradient problem is unique and equals $\alpha_1$ on the curvilinear triangle $p_1 p_2 p_3$ and $0$ and $\alpha_2$ in the respective flops; \\
(3) $h$ is the function from the Brothers example, see \cite[Example 2.7]{MRL}. It is given by the formula

$$h(x,y) = \twopartdef{x^2 - y^2 + 1}{|x| > \frac{1}{\sqrt{2}}}{x^2 - y^2 - 1}{|x| < \frac{1}{\sqrt{2}}.} $$
Then $u \in BV(\Omega)$ is a function of least gradient if and only if 

$$\threepartdef{2x^2}{|x| > \frac{1}{\sqrt{2}}}{\lambda}{|x|, |y| < \frac{1}{\sqrt{2}}}{-2y^2}{|y| > \frac{1}{\sqrt{2}},}$$
where $\lambda \in [-1,1]$.
\end{prz}

A new type of example is the one presented in Section 4.1. There, we witness the phenomenon of breaking of a level set into multiple parts. Of course it can be reversed, i.e. take $u_0$ to be the function which takes two values on the hexagon $H$ and $u$ the function which takes one value; in that case the two level sets of $u_0$ merge into a single level set of $u$. Finally, the following example considers a three-dimensional setting with axial symmetry.

\begin{prz}
Let $\Omega = B(0,1) \subset \mathbb{R}^3$. Take the boundary data to be

$$h(x, y, z) = \twopartdef{1}{|z| > a}{-1}{|z| < a,}$$
where the constant $a$ is chosen so that the two circles which are intersections of $\Omega$ and the planes $\{ z = \pm a \}$ have the same area as the catenoid spanned by them. Using Proposition \ref{stw:jakwygladaroz} and the axial symmetry which helps us solve the Plateau problem we prove that

$$u = \threepartdef{1}{|z| > a}{\lambda}{|z| < a, \text{inside the catenoid}}{-1}{|z| < a, \text{outside the catenoid},}$$
where $\lambda \in [-1,1]$. Moreover, we may take a closer look at the interlacing condition: the set $\{ u \geq t \}$ is the catenoid and the set $\{ u \leq t \}$ is the two circles. The two circles do not intersect and the catenoid intersects the circles at the boundary, so the interlacing condition is satisfied.
\end{prz}


\section{Selection criterion for minimizers}

The strain-gradient plasticity model, as introduced in \cite{ACGR}, is a problem of minimalization of a functional

$$ \widetilde{F_1}(u) = \int_{\Omega} (u^2 + |\nabla u|^2)^{\frac{1}{2}} dx,$$
well-defined over $W^{1,1}(\Omega)$. In the literature, for example see \cite{ACGR}, this functional is minimized with respect to two contraints: the Dirichlet boundary conditions and a condition on the total mass of the solution.

Here we want to introduce a parameter $\varepsilon$ and examine the behavior of minimizers for small $\varepsilon$. For Dirichlet boundary data we define a functional $\widetilde{F_\varepsilon}$ over $L^1(\Omega)$

$$ \widetilde{F_\varepsilon}(u) = \twopartdef{\int_{\Omega} (\varepsilon u^2 + |Du|^2)^{\frac{1}{2}} dx}{u \in W^{1,1}(\Omega), \quad Tu = f}{+ \infty}{\mathrm{otherwise}.}$$

As it turns out even for the simplest possible boundary data, this functional may have no minimizers in $L^1(\Omega)$. We may derive its lower semicontinous envelope similarly as it was calculated in \cite[Section 7]{ACGR} for $\varepsilon = 1$:

$$ F_\varepsilon(u) = \twopartdef{\int_{\Omega} (\varepsilon u^2 + |\nabla u|^2)^{\frac{1}{2}} dx + \int_\Omega |D^s u| + \int_{\partial\Omega} |Tu - f|}{u \in BV(\Omega)}{+ \infty}{\mathrm{otherwise}.}$$

Here we focus on the relationship between this functional and the functional $F$, the relaxed functional in the least gradient problem, namely

$$ F(u) = \twopartdef{\int_{\Omega} |D u| + \int_{\partial\Omega} |Tu - f|}{u \in BV(\Omega)}{+ \infty}{\mathrm{otherwise}.}$$

Firstly, we prove $\Gamma-$convergence of $F_\varepsilon$ (and a similar functional $G_{p, \varepsilon}$) to $F$ and some of its consequences. Secondly, we shall see that minimizers of $G_{p, \varepsilon}$ converge in $L^p$ to minimizers of $F$ which have the smallest norm in $L^p$; this provides a selection criterion for least gradient functions with prescribed boundary conditions, as in general the solutions for Dirichlet least gradient problem may be not unique. Finally, we shall discuss some stronger modes of convergence of these minimizers. We start with recalling the notion of $\Gamma-$convergence:

\begin{dfn}
Let $F, F_n: X \rightarrow [0, \infty]$ be a sequence of functionals on a topological space $X$. We say that the sequence $F_n$ $\Gamma-$converges to $F$, what we denote by $\Gamma-\lim_{n \rightarrow \infty} F_n = F$, if the following two conditions are satisfied: \\
(1) For every sequence $x_n \in X$ such that $x_n \rightarrow x$ in $X$ we have 
$$F(x) \leq \liminf_{n \rightarrow \infty} F_n(x_n);$$
(2) For every $x \in X$ there exists a sequence $x_n \rightarrow x$ in $X$ such that
$$F(x) \geq \limsup_{n \rightarrow \infty} F_n(x_n).$$
We extend this notion for continuous families of parameters in the obvious way: $F_\varepsilon$ $\Gamma-$converges to $F$ as $\varepsilon \rightarrow 0$, if it $\Gamma-$converges for every subsequence. Furthermore, cluster points of minimizers of $F_n$ are minimizers of $F$.
\end{dfn}

\begin{stw}
$\Gamma-\lim_{\varepsilon \rightarrow 0} F_\varepsilon = F$.
\end{stw}

\begin{dd}
We have to check the two conditions in the definition of $\Gamma-$convergence.

(1) We show that for any sequence $u_n \rightarrow u$ in $L^1(\Omega)$ and any sequence $\varepsilon_n \rightarrow 0$ we have $F(u) \leq \liminf_{n \rightarrow \infty} F_{\varepsilon_n}(u_n)$.

$$ \liminf_{n \rightarrow \infty} F_{\varepsilon_n}(u_n) =  \liminf_{n \rightarrow \infty} \int_{\Omega} (\varepsilon_n u_n^2 + |\nabla u_n|^2)^{\frac{1}{2}} dx + \int_\Omega |D^s u_n| + \int_{\partial\Omega} |Tu_n - f| \geq $$
$$  \geq \liminf_{n \rightarrow \infty} \int_{\Omega} |\nabla u_n| dx + \int_\Omega |D^s u_n| + \int_{\partial\Omega} |Tu_n - f| =  \liminf_{n \rightarrow \infty} F(u_n) \geq F(u).$$
The first inequality follows from a pointwise inequality between functions under the integral. The second inequality follows from lower semicontinuity of $F$.

(2) We show that for any function $u \in L^1(\Omega)$ and any sequence $\varepsilon_n \rightarrow 0$ there exists a sequence $u_n \rightarrow u$ such that $F(u) \geq \limsup_{n \rightarrow \infty} F_{\varepsilon_n}(u_n)$.

If $u \notin BV(\Omega)$, the inequality is obvious. If $u \in BV(\Omega)$, take any sequence $u_n$ converging strictly to $u$, i.e. $u_n \rightarrow u$ in $L^1(\Omega)$ and $\int_\Omega |Du_n| \rightarrow \int_\Omega |Du|$. In particular, $\int_\Omega |u_n| dx \leq M$. Then

$$ \limsup_{n \rightarrow \infty} F_{\varepsilon_n}(u_n) =  \limsup_{n \rightarrow \infty} \int_{\Omega} (\varepsilon_n u_n^2 + |\nabla u_n|^2)^{\frac{1}{2}} dx + \int_\Omega |D^s u_n| + \int_{\partial\Omega} |Tu_n - f| \leq $$
$$ \leq \limsup_{n \rightarrow \infty} \int_{\Omega} (\sqrt{\varepsilon_n} |u_n| + |\nabla u_n|) dx + \int_\Omega |D^s u_n| + \int_{\partial\Omega} |Tu_n - f| \leq $$
$$ \leq \limsup_{n \rightarrow \infty} \sqrt{\varepsilon_n} M dx + \int_\Omega |D u_n| + \int_{\partial\Omega} |Tu_n - f| = 0 + \limsup_{n \rightarrow \infty} F(u_n) = F(u).$$
The first inequality follows from a pointwise inequality between functions under the integral. The second inequality follows the upper bound on $L^1$ norms of $u_n$. The limit of $F(u_n)$ equals $F(u)$ because of strict convergence and continuity of trace in the strict topology. \qed
\end{dd}

\begin{uw}
Note that in particular we proved that for strict convergence $u_n \rightarrow u$ we have $F_{\varepsilon_n}(u_n) \rightarrow F(u)$.
\end{uw}

From $\Gamma-$convergence of $F_\varepsilon$ to $F$ it follows that if $u_n$ is a minimizer of $F_{\varepsilon_n}$, then every cluster point of the sequence $u_n$ is a minimizer of $F$. We shall see that we have a common bound in $BV$ norm for minimizers of $F_{\varepsilon}$ for $\varepsilon \leq 1$, so there is a convergent subsequence in $L^1(\Omega)$. 

\begin{stw}\label{zbieznoscl1}
Let $u_n$ be a sequence of minimizers of $F_{\varepsilon_n}$, $\varepsilon_n \rightarrow 0$. We may assume that $\varepsilon \leq 1$. Then there is a convergent subsequence $u_{n_k} \rightarrow u$ in $L^1(\Omega)$.
\end{stw}
 
\begin{dd}
Notice that for $f \in L^1(\partial\Omega)$ we have $F_1(v \equiv 0) = \int_\Omega 0 + \int_{\partial\Omega} |f| < \infty$. Then
$$ \int_\Omega |Du_n| \leq F(u_n) \leq F_{\varepsilon_n}(u_n) \leq F_{\varepsilon_n}(v \equiv 0) \leq F_1(v \equiv 0) < \infty. $$
Thus the total variations of $u_n$ are uniformly bounded. Together with the Dirichlet boundary condition it implies a common bound in $L^1$ norm, so also in $BV$ norm: take an extension of $u_n$ on some ball $B(0,R)$ such that $\Omega \subset \subset B(0,R)$ defined by the formula
$$\widetilde{u_n}(x) = \twopartdef{u_n(x)}{x \in \Omega}{0}{\mathrm{otherwise}.}$$
We apply the Poincar\'{e} inequality to $\widetilde{u_n}$ (note that $\widetilde{u_n}$ has compact support). Thus

$$ \int_\Omega |u_n| = \int_{B(0,R)} |\widetilde{u_n}| \leq C |D\widetilde{u_n}|(B(0,R)) = C |Du_n|(\Omega) + 0 + C \int_{\partial\Omega} |Tu_n| \leq $$
$$ \leq C |Du_n|(\Omega) + C \int_{\partial\Omega} |Tu_n - f| + C \int_{\partial_\Omega} |f| = C F(u_n) + C \int_{\partial_\Omega} |f| \leq$$
$$ \leq C F_1(v \equiv 0) + C F_1(v \equiv 0) = 2C F_1(v \equiv 0) < \infty.$$
It follows that $\|u_n \|_{BV} \leq (2C + 1) F_1(v \equiv 0) < \infty$, so it has a convergent subsequence $u_{n_k} \rightarrow u$ in $L^1(\Omega)$. \qed
\end{dd}

The following result shows that the convergence guaranteed by Proposition \ref{zbieznoscl1} is sometimes in fact not only in $L^1(\Omega)$, but in strict topology of $BV(\Omega)$.

\begin{stw}\label{stw:zbieznoscscisla1}
Let $u_n$ be a sequence of minimizers of $F_{\varepsilon_n}$, $\varepsilon_n \rightarrow 0$. Let $u_n \rightarrow u$ in $L^1(\Omega)$; in particular $u$ is a minimizer of $F$. Then: \\
(1) $F(u_n) \rightarrow F(u)$; \\
(2) If $Tu = f$, then $u_n \rightarrow u$ in the strict topology of $BV(\Omega)$.
\end{stw}

\begin{dd}
(1) Because $u_n$ are minimizers of $F_{\varepsilon_n}$ and $u$ is a minimizer of $F$, we have

$$F(u) \leq F(u_n) \leq F_{\varepsilon_n}(u_n) \leq F_{\varepsilon_n}(u) \rightarrow F(u).$$
(2) Because $u_n$ are minimizers of $F_{\varepsilon_n}$, we have

$$ \int_\Omega |Du_n| \leq F_{\varepsilon_n}(u_n) \leq F_{\varepsilon_n}(u) \leq \int_\Omega \sqrt{\varepsilon_n} |u| + \int_\Omega |Du| + \int_{\partial\Omega} 0,$$
so

$$ \limsup_{n \rightarrow \infty} \int_\Omega |Du_n| \leq \int_\Omega |Du|.$$
By lower semicontinuity of the total variation we obtain the opposite inequality, so $\lim_{n \rightarrow \infty} \int_\Omega |Du_n| = \int_\Omega |Du|.$ \qed
\end{dd}


However, looking at the functional $F_{\varepsilon}$ gives us little information about pointwise properties of the approximating sequence $u_n$. It also gives us convergence to some minimizer of $F$, while we want our sequence to choose one particular element of $\arg\min F$. To this end, let us define for $1 \leq p < \frac{k}{k-1}$ an auxiliary functional $G_{p,\varepsilon}$:

$$ G_{p,\varepsilon}(u) = \twopartdef{(\int_{\Omega} \sqrt{\varepsilon} |u|^p )^{\frac{1}{p}} + \int_\Omega |D u| + \int_{\partial\Omega} |Tu - f|}{u \in BV(\Omega)}{+ \infty}{\mathrm{otherwise}.}$$
In other words, we have $G_{p,\varepsilon}(u) = \sqrt[2p]{\varepsilon} \| u \|_p + F(u)$. Using the continuous embedding of $BV(\Omega)$ into $L^p(\Omega)$, we see that all the above results hold also for $G_{p,\varepsilon}$ with an analogous proof. Now we shall see that $G_{p,\varepsilon}$ provides a selection criterion for minimizers of $F$:

\begin{tw}\label{najmniejszanorma}
Let $v_n \in \arg\min G_{p, \varepsilon_n}$ and $\varepsilon_n \rightarrow 0$. Suppose that $v_n \rightarrow v$ in $L^p(\Omega)$. By $\Gamma-$convergence of $G_{\varepsilon_n}$ we have $v\in \arg\min F$. Then $v$ is an element with the smallest $L^p$ norm among minimizers of $F$.
\end{tw}

\begin{dd}
Suppose that $u$ is another minimizer of $F$, which has smaller $L^p$ norm than $v$. Let $\delta < \frac{\|v \|_p - \| u \|_p}{2}$. Fix $n$ big enough, namely let $|\| v_n \|_p - \| v \|_p| < \delta$. As $v_n$ are minimizers of $G_{\varepsilon_n}$, we have

$$ 0 \geq G_{p, \varepsilon_n}(v_n) - G_{p, \varepsilon_n}(u) = \sqrt[2p]{\varepsilon} \| v_n \|_p + F(v_n) - \sqrt[2p]{\varepsilon} \| u \|_p - F(u) \geq$$
$$ \geq \sqrt[2p]{\varepsilon} \| v_n \|_p - \sqrt[2p]{\varepsilon} \| u \|_p  = \sqrt[2p]{\varepsilon} (\| v_n \|_p - \|v \|_p) + \sqrt[2p]{\varepsilon} (\|v \|_p - \| u \|_p ) \geq - \sqrt[2p]{\varepsilon} \delta + 2 \sqrt[2p]{\varepsilon} \delta > 0,$$
contradiction. Thus $v_n$ cannot converge to an element which does not have smallest $L^p$ norm. \qed
\end{dd}

Let us note that compact embedding of $BV(\Omega)$ into $L^p(\Omega)$ implies that the sequence $u_n$, due to its boundedness in $BV(\Omega)$, is convergent in $L^p(\Omega)$ on some subsequence. As the natural underlying space for $BV(\Omega)$ is $L^1(\Omega)$, it is tempting to consider only $p = 1$; however, we do not know if the minimizer of $F$ with the smallest norm in $L^1$ is unique, while for $p > 1$ it is unique (see later in Proposition \ref{smallestpnorm}). Furthermore, $G_{p, \varepsilon}$ have unique minimizers for $p > 1$, as they are strictly convex; it does not apply to $p = 1$.

However, convergence in $L^1$ is quite weak, so a natural question is if some stronger mode of convergence might be at play. The natural candidate is strict convergence; however, for boundary data with a constant sign we may prove a much stronger result.

\begin{stw}\label{stw:monotone}
Let $f \in L^1(\partial\Omega)$ be nonnegative. Let $\varepsilon_1 > \varepsilon_2$. Then any minimizer of $G_{p, \varepsilon_1}$ is pointwise smaller than any minimizer of $G_{p, \varepsilon_2}$, i.e. let $u_1 \in \arg\min G_{p, \varepsilon_1}$ and $u_2 \in \arg\min G_{p, \varepsilon_2}$. Then $u_2 \geq u_1$. 
\end{stw}

\begin{dd}
In the beginning, let us note that as $f$ is nonnegative, $u_1$ and $u_2$ are as well: it is enough to compare the value of $G_{p, \varepsilon_i}$ on $u_i$ and $\max(u_i,0)$.

Our starting point is the inequality

\begin{equation}\label{eq:porownanieGepsilon}
G_{p, \varepsilon_1}(u_1) + G_{p, \varepsilon_2}(u_2) \leq G_{p, \varepsilon_1}(\min(u_1,u_2))  + G_{p, \varepsilon_2}(\max(u_1, u_2)),
\end{equation}
which is automatically fulfilled, as $u_1$ and $u_2$ are minimizers of $G_{p, \varepsilon_1}$ and $G_{p, \varepsilon_2}$ respectively. Our goal is to prove the opposite inequality and under what conditions is it strict.

Now we expand the left hand side of the above inequality:

$$G_{p, \varepsilon_1}(u_1) + G_{p, \varepsilon_2}(u_2) = \sqrt[2p]{\varepsilon_1} \|u_1 \|_p + \int_\Omega |Du_1| + $$
$$+ \int_{\partial\Omega} |Tu_1 - f| + \sqrt[2p]{\varepsilon_2} \|u_2 \|_p + \int_\Omega |Du_2| + \int_{\partial\Omega} |Tu_2 - f|$$
and the right hand side:

$$G_{p, \varepsilon_1}(\min(u_1, u_2)) + G_{p, \varepsilon_2}(\max(u_1, u_2)) = \sqrt[2p]{\varepsilon_1} \|\min(u_1, u_2) \|_p + $$
$$+\int_\Omega |D\min(u_1, u_2)| + \int_{\partial\Omega} |T\min(u_1, u_2) - f| + \sqrt[2p]{\varepsilon_2} \|\max(u_1, u_2) \|_p + $$
$$+\int_\Omega |D\max(u_1, u_2)| + \int_{\partial\Omega} |T\max(u_1, u_2) - f|.$$
Firstly, let us recall that Lemma \ref{min+max} states that

$$ \int_\Omega |D\max(u,v)| + \int_\Omega |D\min(u,v)| \leq \int_\Omega |Du| + \int_\Omega |Dv|.$$
Secondly, see that Lemma \ref{traceofminimum} implies that

$$ \int_{\partial\Omega} |T\max(u_1, u_2) - f| + \int_{\partial\Omega} |T\min(u_1, u_2) - f| = \int_{\partial\Omega} |T u_1 - f| + \int_{\partial\Omega} |T u_2 - f|$$
as we have pointwise equality $\mathcal{H}^{k-1}$-a.e. Thus most of summands in (\ref{eq:porownanieGepsilon}) cancel out and it reduces to the following inequality:

$$\sqrt[2p]{\varepsilon_1} \|u_1 \|_p + \sqrt[2p]{\varepsilon_2} \|u_2 \|_p \leq \sqrt[2p]{\varepsilon_1} \|\min(u_1, u_2) \|_p + \sqrt[2p]{\varepsilon_2} \|\max(u_1, u_2) \|_p.$$
As $u_1, u_2$ are nonnegative, we may expand the left hand side in the following way:

$$\sqrt[2p]{\varepsilon_1} \|u_1 \|_p + \sqrt[2p]{\varepsilon_2} \|u_2 \|_p = \sqrt[2p]{\varepsilon_1} \int_0^\infty p t^{p-1} |\{ u_1 > t \}| dt + $$
$$+ \sqrt[2p]{\varepsilon_2} \int_0^\infty p t^{p-1} |\{ u_2 > t \}| dt = \sqrt[2p]{\varepsilon_1} \int_0^\infty p t^{p-1} |\{ u_1 > t \} \cap \{ u_2 > t \}| dt +$$
$$ + \sqrt[2p]{\varepsilon_1} \int_0^\infty p t^{p-1} |\{ u_1 > t \} \backslash \{ u_2 > t \}| dt + \sqrt[2p]{\varepsilon_2} \int_0^\infty p t^{p-1} |\{ u_2 > t \} \backslash \{ u_1 > t \}| dt +$$ $$+ \sqrt[2p]{\varepsilon_2} \int_0^\infty p t^{p-1} |\{ u_1 > t \} \cap \{ u_2 > t \}| dt.$$
And the right hand side in the following way:

$$\sqrt[2p]{\varepsilon_1} \|\min(u_1, u_2) \|_p + \sqrt[2p]{\varepsilon_2} \|\max(u_1, u_2) \|_p = \sqrt[2p]{\varepsilon_1} \int_0^\infty p t^{p-1} |\{ \min(u_1, u_2) > t \}| dt + $$
$$+ \sqrt[2p]{\varepsilon_2} \int_0^\infty p t^{p-1} |\{ \max(u_1, u_2) > t \}| dt = \sqrt[2p]{\varepsilon_1} \int_0^\infty p t^{p-1} |\{ u_1 > t \} \cap \{ u_2 > t \}| dt +$$
$$ + \sqrt[2p]{\varepsilon_2} \int_0^\infty p t^{p-1} |\{ u_1 > t \} \backslash \{ u_2 > t \}| dt + \sqrt[2p]{\varepsilon_2} \int_0^\infty p t^{p-1} |\{ u_2 > t \} \backslash \{ u_1 > t \}| dt +$$ $$+ \sqrt[2p]{\varepsilon_2} \int_0^\infty p t^{p-1} |\{ u_1 > t \} \cap \{ u_2 > t \}| dt.$$
Again, most of the summands cancel out and we are left with

$$(\sqrt[2p]{\varepsilon_2} - \sqrt[2p]{\varepsilon_1}) \int_0^\infty p t^{p-1} |\{ u_1 > t \} \backslash \{ u_2 > t \}| dt \geq 0,$$
which implies that for almost every $t$ the Lebesgue measure of the set $\{ u_1 > t \} \backslash \{ u_2 > t \}$ is zero, so $u_2 \geq u_1$ a.e. \qed 
\end{dd}

At this point, let us clearly state a few implications of the above result.

\begin{wn}
In particular, if $\varepsilon_n$ goes monotonically to zero, then for nonnegative boundary data every sequence of minimizers of $G_{p, \varepsilon_n}$ is convergent to $u$ without the need of choosing a subsequence. Furthermore, let us look closer at the inequality $G_{p, \varepsilon_2}(u_2) \leq G_{p, \varepsilon_2}(u_1)$ (true by definition of $u_2$). After expanding both sides we get

$$\sqrt[2p]{\varepsilon_2} \|u_2 \|_p + F(u_2) \leq \sqrt[2p]{\varepsilon_2}  \|u_1 \|_p + F(u_1),$$
so in view of Proposition \ref{stw:monotone} this implies that $F(u_1) \geq F(u_2)$. Thus the sequence $F(u_n)$ is decreasing. By Proposition \ref{stw:zbieznoscscisla1} it converges to $F(u)$.
\end{wn}

The fact that by Proposition \ref{stw:monotone} $u_n$ is an increasing sequence allows us to prove an improved version of Proposition \ref{stw:zbieznoscscisla1}. 

\begin{wn}
Take an increasing sequence $u_n \rightarrow u$ in $L^p(\Omega)$ as mentioned in the previous Corollary. Suppose that $Tu \leq f$. Then $u_n \rightarrow u$ in the strict topology of $BV(\Omega)$. 
\end{wn}

\begin{dd}
We proceed similarly to the proof of Proposition \ref{stw:zbieznoscscisla1}:

$$\int_\Omega |Du_n| + \int_{\partial\Omega} |Tu_n - f| = F(u_n) \leq G_{p, \varepsilon_n}(u_n) \leq G_{p, \varepsilon_n}(u) \rightarrow$$
$$\rightarrow F(u) =  \int_\Omega |Du| + \int_{\partial\Omega} |Tu - f|,$$
so

$$\limsup_{n \rightarrow \infty} \int_\Omega |Du_n| + \limsup_{n \rightarrow \infty} \int_{\partial\Omega} |Tu_n - f| \leq \int_\Omega |Du| + \int_{\partial\Omega} |Tu - f|.$$
On the other hand, by monotonicity of $u_n$ we have $Tu_n \leq Tu \leq f$. In particular, $\int_{\partial\Omega} |Tu_n - f| \leq \int_{\partial\Omega} |Tu - f|$. This coupled with the lower semicontinity of the total variation gives us

$$\liminf_{n \rightarrow \infty} \int_\Omega |Du_n| + \liminf_{n \rightarrow \infty} \int_{\partial\Omega} |Tu_n - f| \geq \int_\Omega |Du| + \int_{\partial\Omega} |Tu - f|.$$
This means that every inequality in an equality, in particular $\lim_{n \rightarrow \infty} \int_\Omega |Du_n| = \int_\Omega |Du|$. \qed
\end{dd}


As it was mentioned above, we are going to take advantage of the fact that for $1 < p < \frac{k}{k-1}$ there is a unique minimizer of $F$ in $L^p(\Omega)$. This will give us convergence on the whole sequence of minimizers of $G_{p, \varepsilon_n}$. Moreover, it turns out that Theorem \ref{uniqueness} helps us to estabilish a similar claim for minimizers of $F$ which attain the trace $f$ also for $p = 1$.

\begin{stw}\label{smallestpnorm}
Let $X$ be the set of minimizers of $F$. Then $X$ is a compact convex set in $L^p(\Omega)$, where $1 \leq p < \frac{k}{k-1}$. In particular it has a unique element of the smallest $p-$norm for $1 < p < \frac{k}{k-1}$.
\end{stw}

\begin{dd}
As $F$ is convex, the arithmetic mean of minimizers is also a minimizer, so $X$ is convex. As $F$ is lower semicontinuous, the set of minimizers is closed in $L^1(\Omega)$ (as the limit of minimizers attains the same value of $F$), so it is closed in $BV(\Omega)$ and by continuity of the embedding into $L^p(\Omega)$ for $1 \leq p \leq \frac{k}{k-1}$ it is closed in $L^p(\Omega)$ for $1 \leq p \leq \frac{k}{k-1}$. It is a bounded set in every $L^p(\Omega)$ for $1 \leq p \leq \frac{k}{k-1}$, as it is bounded in $BV(\Omega)$: firstly, if $u$ is a minimizer of $F$, then $\int_\Omega |Du| \leq F(u) = m$, the minimal value of $F$ (note that also $\int_{\partial\Omega} |Tu - f| \leq m$). Secondly, let us extend $u$ by $0$ on some ball $B(0,r)$ including $\Omega$. From the Poincar\'{e} inequality we have

$$ \| u \|_1 \leq C( \int_\Omega |Du| + \int_{\partial\Omega} |Tu|) \leq C( \int_\Omega |Du| +  \int_{\partial\Omega} |Tu - f| + \int_{\partial\Omega} |f|) \leq$$
$$ \leq C(2m + \int_{\partial\Omega} |f|) = M,$$
so $X$ is a bounded set in $BV(\Omega)$. Thus $X$ is bounded and closed in $L^p$. For $1 \leq p < \frac{k}{k-1}$ it is compact in $L^p$, so for $1 < p < \frac{N}{N-1}$ it has a unique element of the smallest norm. \qed
\end{dd}

\begin{wn}
Thus for $1 < p < \frac{k}{k-1}$ the minimizers of $G_{p, \varepsilon_n}$, $u_n$, converge to $u$, converge the element of the smallest $p-$norm of $X$ not only on some subsequence, but on the whole sequence: it is a consequence of the fact that in metric spaces (and $BV(\Omega)$ endowed with strict topology is metrizable) if we can from every subsequence $x_{n_l}$ extract a subsubsequence $x_{n_{l_m}} \rightarrow x$, then $x_n \rightarrow x$. It provides a selection criterion for elements of $X$.
\end{wn}

\begin{wn}
In a slightly different case, where $X$ is the set of minimizers of $F$ with trace $f$ (the boundary condition is met in the trace sense), Proposition \ref{smallestpnorm} also holds. This is a consequence of Theorem \ref{uniqueness}.
\end{wn}

\begin{dd}
The proof of convexity, boundedness and compactness does not change. We only have to prove that $X$ is closed in $L^p(\Omega)$ (it is enough to prove closedness in $L^1(\Omega)$). 

Take a sequence $u_n$ of least gradient functions with trace $f$ which converges to $u$ in $L^1(\Omega)$. By Miranda's Theorem, see \cite[Theorem 3]{Mir} $u$ is a function of least gradient. Now, by Theorem \ref{uniqueness} the functions $u_m$ and $u_n$ differ only on some set $C_{nm}$, on which both functions are locally constant. As both functions have the same trace, we have $\mathcal{H}^1(\overline{C_{nm}} \cap \partial\Omega) = 0$. If we denote $C = \bigcup_{n = 1}^{\infty} C_{nm}$, then 
$$\mathcal{H}^{k-1}(\overline{C} \cap \partial\Omega) \leq \sum_{n, m} \mathcal{H}^{k-1}(\overline{C_{nm}} \cap \partial\Omega) = 0.$$
On $\Omega \backslash C$ the sequence $u_n$ is constant, so on this set $u = u_n$. As $\mathcal{H}^{k-1}(\overline{C} \cap \partial\Omega) = 0$, the boundary of $\Omega \backslash C$ is the whole $\partial\Omega$. This means that $Tu = Tu_n = f$, so $X$ is a closed set. \qed
\end{dd}

Let us conclude this section with noticing that Theorem \ref{uniqueness} together with the analysis in Section 4 implies that the element of $X$ with the smallest $1-$norm is the same as the element of $X$ with the smallest $p-$norm, in particular it is unique. Thus the functional $G_{p,\varepsilon}$ produces a selection criterion for elements of $X$ also for $p=1$.

\textbf{Acknowledgement.} I would like to thank my PhD advisor, Piotr Rybka, for numerous discussions concerning this paper, and Lorenzo Giacomelli for suggesting the possible connection between the least gradient problem and the strain-gradient plasticity model.

\end{document}